\documentclass[leqno,10pt]{amsart}
\usepackage{amsmath,amscd,amsthm,amsxtra}
\usepackage{epsfig,graphics,color,colortbl}
\usepackage{amssymb,latexsym}
\usepackage{mathabx}
\usepackage{mathrsfs,eucal,upgreek}
\usepackage[poly,all]{xy}
\usepackage{hyperref}
\usepackage{tikz-cd}
\usepackage{arydshln}
\usepackage{ytableau}
\usepackage{adjustbox}
\usepackage{listings}
\usepackage{dirtytalk}
\usepackage[draft]{todonotes}
\usepackage{datetime}

\usepackage{marginnote}
\usepackage[normalem]{ulem}

\newtheorem{thm}{\bf Theorem}[section]
\newtheorem{df}[thm]{\bf Definition}
\newtheorem{prop}[thm]{\bf Proposition}
\newtheorem{cor}[thm]{\bf Corollary}
\newtheorem{lem}[thm]{\bf Lemma}
\newtheorem{rem}[thm]{\bf Remark}
\newtheorem{ex}[thm]{\bf Example}

\numberwithin{equation}{section}

\makeatletter
\renewcommand{\tocsection}[3]{%
  \indentlabel{\@ifnotempty{#2}{\bfseries\ignorespaces#1 #2.\quad}}\bfseries#3}
\renewcommand{\tocsubsection}[3]{%
  \indentlabel{\@ifnotempty{#2}{\ignorespaces#1 #2.\quad}}#3}
\newcommand\@dotsep{4.5}
\def\@tocline#1#2#3#4#5#6#7{\relax
  \ifnum #1>\c@tocdepth % then omit
  \else
    \par \addpenalty\@secpenalty\addvspace{#2}%
    \begingroup \hyphenpenalty\@M
    \@ifempty{#4}{%
      \@tempdima\csname r@tocindent\number#1\endcsname\relax
    }{%
      \@tempdima#4\relax
    }%
    \parindent\z@ \leftskip#3\relax \advance\leftskip\@tempdima\relax
    \rightskip\@pnumwidth plus1em \parfillskip-\@pnumwidth
    #5\leavevmode\hskip-\@tempdima{#6}\nobreak
    \leaders\hbox{$\m@th\mkern \@dotsep mu\hbox{.}\mkern \@dotsep mu$}\hfill
    \nobreak
    \hbox to\@pnumwidth{\@tocpagenum{\ifnum#1=1\bfseries\fi#7}}\par% <-- \bfseries for \section page
    \nobreak
    \endgroup
  \fi}
\AtBeginDocument{%
\expandafter\renewcommand\csname r@tocindent0\endcsname{0pt}
}
\def\l@subsection{\@tocline{2}{0pt}{2.5pc}{5pc}{}}
\makeatother
\newcommand{\B}{\mathbf{B}}

\newcommand{\pf}{\noindent{\bfseries Proof. }}
\newcommand{\ov}{\overline}

\newcommand{\Z}{\mathbb{Z}}
\newcommand{\C}{\mathbb{C}}

\newcommand{\te}{\widetilde{e}}

\newcommand{\g}{\mathfrak{g}}

\newcommand{\mc}{\mathcal}
\newcommand{\mf}{\mathfrak}
\newcommand{\La}{\Lambda}

\newcommand{\ep}{\epsilon}
\newcommand{\vep}{\varepsilon}

\newcommand*\oline[1]{%
   \vbox{%
     \hrule height 0.85pt%                  % Line above with certain width
     \kern0.25ex%                          % Distance between line and content
     \hbox{%
       \kern -0.05em%                        % Distance between content and left side of box, negative values for lines shorter than content
       \ifmmode#1\else\ensuremath{#1}\fi%  % The content, typeset in dependence of mode
       \kern -0em%                        % Distance between content and left side of box, negative values for lines shorter than content
     }% end of hbox
   }% end of vbox
}

\newcommand{\al}{{\alpha}}

\newcommand{\xz}{\mathsf{x}_0}

\newcommand{\bR}{\mathbf{k}}

\newcommand{\pfe}{\mathsf{e}}
\newcommand{\pfk}{\mathsf{k}}

\newcommand{\E}{\mathsf{E}} % real root vector E
\newcommand{\F}{\mathsf{F}} % real root vector F
\newcommand{\V}{{\bf 1}} % (loop) highest weight vector

\newcommand{\ot}{\otimes}
\newcommand{\wtd}{\widetilde}

\newenvironment{redm}{\relax\color{red}}{\relax}
\newenvironment{bluem}{\relax\color{blue}}{\relax}
\newenvironment{jaunem}{\relax\color{pink}}{\relax}

\newcommand{\beb}{\begin{bluem}}
	\newcommand{\eb}{\end{bluem}}

\newcommand{\ber}{\begin{redm}}
\newcommand{\er}{\end{redm}}

\newcommand{\bej}{\begin{jaunem}}
\newcommand{\ej}{\end{jaunem}}

\newcommand{\berE}[1]{\begin{redm}{}\marginnote{\fbox{\scshape\lowercase{E}}}	#1}  % Euiyong
\newcommand{\bejE}[1]{\begin{jaunem}{}\marginnote{\fbox{\scshape\lowercase{E}}}	#1}  % Euiyong
\newcommand{\berJ}[1]{\begin{redm}{}\marginnote{\fbox{\scshape\lowercase{J}}}	#1}  % Il-Seung
\newcommand{\bejJ}[1]{\begin{jaunem}{}\marginnote{\fbox{\scshape\lowercase{J}}}	#1}  % Il-Seung

\begin{document}

\title[Prefundamental representations for types $A_n^{(1)}$ and $D_n^{(1)}$]
{Unipotent quantum coordinate ring and prefundamental representations for types $A_n^{(1)}$ and $D_n^{(1)}$}

\author{IL-SEUNG JANG}

\address{Department of Mathematical Sciences, Seoul National University, Seoul 08826, Republic of Korea}
\email{is\_jang@snu.ac.kr}

\author{JAE-HOON KWON}

\address{Department of Mathematical Sciences and RIM, Seoul National University, Seoul 08826, Republic of Korea}
\email{jaehoonkw@snu.ac.kr}

\author{EUIYONG PARK}
\address{Department of Mathematics, University of Seoul, Seoul 02504, Republic of Korea}
\email{epark@uos.ac.kr}

\keywords{quantum affine algebra, prefundamental representation, unipotent quantum coordinate ring}
\subjclass[2010]{17B37, 22E46, 05E10}

\thanks{I.-S. Jang and J.-H. Kwon are supported by the National Research Foundation of Korea(NRF) grant funded by the Korea government(MSIT) (No.2019R1A2C1084833 and No. 2020R1A5A1016126)}
\thanks{
E. Park is supported by the National Research Foundation of Korea (NRF) Grant funded by the Korea Government (MSIT) (NRF- 2020R1F1A1A01065992)}

\begin{abstract}
We give a new realization of the prefundamental representations $L^\pm_{r,a}$ introduced by Hernandez and Jimbo, when the quantum loop algebra $U_q(\mf{g})$ is of types $A_n^{(1)}$ and $D_n^{(1)}$, and the $r$-th fundamental weight $\varpi_r$ for types $A_n$ and $D_n$ is minuscule.
We define an action of the Borel subalgebra $U_q(\mf{b})$ of $U_q(\mf{g})$ on the unipotent quantum coordinate ring associated to the translation by $-\varpi_r$, and show that it is isomorphic to $L^\pm_{r,a}$.
We then give a combinatorial realization of $L^+_{r,a}$ in terms of the Lusztig data of the dual PBW vectors.
\end{abstract}

\maketitle
\setcounter{tocdepth}{2}

\section{Introduction} 
Let $U_q(\mf{g})$ be the quantum loop algebra associated to an affine Kac-Moody algebra $\mf{g}$ of untwisted type.  
Let $I=\{\,1,\dots,n\,\}$ be the index set for the Dynkin diagram of the finite-dimensional simple subalgebra $\mathring{\mf{g}}$ of $\mf{g}$.
Let $x_{i,r}^\pm$, $k^{\pm 1}_{i}$, $h_{i,s}$ ($i\in I, r,s \in\Z, s\neq 0$) denote the Drinfeld generators of $U_q(\mf{g})$ \cite{D2}.

A finite-dimensional irreducible $U_q(\mf{g})$-module (of type I) is generated by an eigenvector $v$ with respect to the currents $(\psi^\pm_i(z))_{i\in I}$ (see \eqref{eq:psi generators}) such that $x^+_{i,r}v=0$ for all $i\in I$, $r\in \Z$.
The eigenvalue $\Psi_i(z)$ of $v$ with respect to $\psi_i^\pm(z)$ is given by a rational functions in $z$ of the following form
\begin{equation*}
\Psi_i(z)= q_i^{{\rm deg}(P_i)}\frac{P_i(q_i^{-1}z)}{P_i(q_iz)},
\end{equation*}
where $P_i(z)$ is a polynomial in $z$ such that $P_i(0)=1$ and $q_i$ is a power of $q$ depending on $i$.
The sequence of  polynomials $P=(P_i(z))_{i\in I}$ parametrizes the finite-dimensional irreducible $U_q(\mf{g})$-modules of type I \cite{CP94, CP}.

Let $U_q(\mf{b})$ be the Borel subalgebra $\mf{b}$ of $\mf{g}$.
In \cite{HJ}, Hernandez and Jimbo introduced a category $\mc{O}$ of $U_q(\mf{b})$-modules, which contains the finite-dimensional $U_q(\mf{g})$-modules. 
An irreducible module in $\mc{O}$ is in general infinite-dimensional, and it is generated by an eigenvector vector $v$ with respect to $(\psi^+_i(z))_{i\in I}$ such that $e_iv=0$ for all $i\in I$ and the eigenvalue $\Psi_i(z)$ for $\psi^+_i(z)$ is given by an arbitrary rational function in $z$. Here $e_i$ is the Chevalley generator of the positive part of $U_q(\mf{g})$.

For a non-zero $a\in \C(q)$ and $r\in I$, the prefundamental representations $L_{r,a}^\pm$ are the $U_q(\mf{b})$-modules corresponding to $(\Psi_i(z))_{i\in I}$ such that
\begin{equation*}
\Psi_r(z) = (1-az)^{\pm 1},
\end{equation*}
and $\Psi_i(z)=1$ elsewhere. 
It is shown in \cite{HJ} that $L_{r, 1}^\pm$ is irreducible and $L_{r, 1}^-$ is a limit of Kirillov-Reshetikhin (simply KR) modules $W^{(r)}_{s,\, q_i^{-2s+1}}$ $(s \geq 1)$, which gives a representation theoretic explanation on the fact \cite{H06,Na03} that the normalized $q$-character of $W^{(r)}_{s,\, q_i^{-2s+1}}$ \cite{FR} has a well-defined limit as $s\rightarrow \infty$. This implies that $L_{r, a}^{\pm}$ are $U_q({\mf b})$-modules in $\mc{O}$, and any irreducible module in $\mc{O}$ is a subquotient of a tensor product of $L_{r, a}^{\pm}$'s, hence parametrized by the $I$-tuples of rational functions $\Psi_i(z)$ for $i \in I$ up to $1$-dimensional $U_q({\mf b})$-modules with a trivial action of $e_i$.

The category $\mc{O}$ plays an important role in generalizing the Baxter's relation \cite{FH}, where the $q$-character of a finite-dimensional $U_q(\mf{g})$-module $V$ is realized as a relation in the Grothendieck ring of $\mc{O}$ including the class of $V$ tensored by a certain tensor product of $L_{r,a}^+$'s.
Furthermore, it is proved in \cite{HL16} that the Grothendieck ring of a certain monoidal subcategory of $\mc{O}$ given by a restriction on the zeros and poles of $\Psi_i(z)$ for all $i \in I$ has a structure of cluster algebra, where the class of prefundamental representations $L_{r,a}^+$ form an initial seed.
 
Despite the importance of $L_{r, a}^{\pm}$, a realization of $L_{r, a}^{\pm}$ does not seem to be known much in general except for special cases.
The purpose of this paper is to give an explicit realization of prefundamental representations $L_{r, a}^{\pm}$ when $\mf{g}=A_n^{(1)}$, $D_n^{(1)}$ and $r\in I$ is minuscule, that is, the corresponding fundamental representation of $\mathring{\mf{g}}$ is minuscule.

Let $\varpi_r$ denote the $r$-th fundamental weight for $\mathring{\mf{g}}$. 
The translation $t_{\varpi_r}$ by $\varpi_r$ in the extended Weyl group of $\mf{g}$ is given by $t_{\varpi_r}=\tau (w_r)^{-1}$ for some $\tau$ and $w_r$, where $\tau$ is an automorphism of the Dynkin diagram and $w_r$ is the element in the Weyl group of $\mf{g}$.
We consider the unipotent quantum coordinate ring associated to $w_r$, which we denote by $U^-_q(w_r)$.
It is a subalgebra of $U_q^-({\mf g})$, the negative part of $U_q({\mf g})$, and this subalgebra is of special importance since it can be viewed as a $q$-analog of the coordinate ring of the unipotent subgroup associated to $w_r$, and has a quantum cluster algebra structure \cite{GLS13}. As a quantum cluster algebra, it also has a monoidal categorification in terms of the representations of quiver Hecke algebras \cite{KKKO}.

We show that there exist two $U_q(\mf{b})$-module structures $\rho^\pm_{r,a}$
on $U^-_q(w_r)$ such that
\begin{equation} \label{eq:main result}
U^-_q(w_r) \cong L_{r, \pm ac_r}^\pm,
\end{equation}
where $c_r$ is a scalar given in \eqref{eq:c_r}.
Note that we may assume that $U^-_q(w_r)$ is a subalgebra of $U_q^-(\mathring{\mf{g}})\subset U_q^-(\mf{g})$ in our case. The action of $e_i$ on $U^-_q(w_r)$ for $i\in I$ is given by the usual $q$-derivation on $U_q^-(\mf{g})$, while the action of $e_0$ on $U^-_q(w_r)$ is given as a right (resp. left) 
multiplication by the dual PBW vector $\xz$ corresponding to the maximal root (see \eqref{eq:xz} for the definition of $\xz$)
by which $U^-_q(w_r)$ is isomorphic to $L_{r, ac_r}^+$ (resp. $L_{r,\,-ac_r}^-$) as a $U_q({\mf b})$-module. Here $e_0$ is the Chevalley generator corresponding to the vertex $0$ in the affine Dynkin diagram of $\mf{g}$.
We remark there is another construction of prefundamental representations $L^-_{1,a}$ and $L^+_{n,a}$ for $A_n^{(1)}$ in terms of $q$-oscillator representations \cite{BGKNR16, BGKNR17}. Our work can be viewed as a generalization of it.

Furthermore, we give an explicit combinatorial realization of $L_{r, a}^+$ in terms of the Lusztig data of the dual PBW vectors. More precisely, we identify the dual PBW vectors of $U^-_q(w_r)$ with the matrices consisting of the Lusztig data of the dual PBW vectors in a natural way and describe the $U_q({\mf b})$-actions explicitly under this identification.

One of the key ingredients in the proof of \eqref{eq:main result} and the above combinatorial realization is that the $q$-derivations on the dual PBW vectors of $U^-_q(w_r)$ are well-understood by the minuscule representation $V(\varpi_r)$ of $U_q(\mathring{\mf g})$. This enables us to  give a simple description of $U_q({\mf b})$-actions on $U^-_q(w_r)$,
in particular, the action of the current $\psi^+_i(z)$ on a maximal vector $1 \in U^-_q(w_r)$, which plays a crucial role in finding its loop highest weight.

Our realization of $L_{r, a}^\pm$ is partly motivated by the results in \cite{K13, K18}, where the affine $A_n^{(1)}$-crystal structure on $U^-_q(w_r)$ is studied in terms of the Lusztig data of the PBW vectors of type $A_n$. More precisely, it is proved in \cite{K13, K18} that the associated affine crystal denoted by $\B^J$ $(J := I \,/\, \{ r \})$ is isomorphic to the limit of the crystals of KR modules $W^{(r)}_s$ as $s \rightarrow \infty$ (see \cite{JK19} for type $D_n^{(1)}$). In particular, the Kashiwara operator $\widetilde{e}_0$ on $\B^J$ is given simply by increasing the multiplicity of the maximal root vector by $1$ in the Lusztig data. This motivated the construction of $L_{r, a}^\pm$ in \eqref{eq:main result}, and now we may understand $\B^J$ as the crystal of $L_{r, 1}^+$.

The paper is organized as follows. In Section \ref{sec:preliminaries}, we briefly review necessary background on quantum loop algebras and the category $\mc{O}$. In Section \ref{sec:UQCR}, we define a $U_q(\mf{b})$-module structure $\rho_{r,a}^+$ of $U_q^-(w_r)$, which belongs to the category $\mc{O}$ (Theorem \ref{thm:main-1}).
In Section \ref{sec:realization}, we compute the eigenvalues of $1 \in U_q^-(w_r)$ with respect to $(\psi_i^+(z))_{i\in I}$ and show that $U_q^-(w_r) \cong L_{r,ac_r}^+$ in case of $\rho_{r,a}^+$ (Theorem \ref{thm:main-2}). Then we define $\rho^-_{r,a}$ in a similar way, and show that $U_q^-(w_r) \cong L_{r,\,-ac_r}^-$ in this case (Theorem \ref{thm:main-3}). 
In Section \ref{sec:combinatorial realization},  we give an explicit realization description of $L_{r, a}^+$ in terms of the Lusztig data of the dual PBW vectors (Theorem \ref{thm:combinatorial realization for type A} for type $A_n^{(1)}$ and Theorem \ref{thm:combinatorial realization for type D} for type $D_n^{(1)}$).

\vskip 3mm
{\em \bf Acknowledgement.} The authors would like to thank D. Hernandez for his interest and comments, and A. V. Razumov for letting us know the references \cite{BGKNR16, BGKNR17}. They also would like to thank the anonymous referees for helpful comments.
%%%%%%%%%%%%%%%%%%%%%%%%%%%%%%%%%%%%%%%%%%%%%%%%%%%%%%%%%%%%%%%%%%%%%%%%

\section{Preliminaries} \label{sec:preliminaries}
\subsection{Notations}
Let $\Z_+$ denote the set of non-negative integers.
Let $A=(a_{ij})_{0\le i,j\le n}$ be the generalized Cartan matrix of symmetric affine type, and let $\mf{g}$ be the affine Kac-Moody algebra associated to $A$. Put $I=\{\,1,\dots,n\,\}$. 
Let $\mathring{A}=(a_{ij})_{i,j\in I}$ be the Cartan matrix of finite type, and let $\mathring{\mf{g}}$ denote the associated finite-dimensional simple Lie algebra.

Let $\{\,\alpha_i\,|\, 0 \le i \le n \,\}$ be the set of simple roots for ${\mf g}$, and let $Q=\bigoplus_{0 \le i \le n}\Z\alpha_i$ be the root lattice. Put $Q_+=\bigoplus_{0\le i\le n}\Z_+\alpha_i$.
Let $ \{\, \Lambda_i \,\mid\, 0 \le i \le n \,\}$ be the fundamental weights for $\g$, and let $P = \bigoplus_{0 \le i \le n} \mathbb{Z}\Lambda_i \bigoplus \mathbb{Z}\delta$ be the weight lattice of $\g$, where $\delta$ is the imaginary null root.
Let us take a nondegenerate symmetric bilinear form $(\ ,\ )$ on $P$
so that $(\alpha_i,\alpha_i)\in 2\Z_{>0}$, $\frac{2(\alpha_i,\alpha_j)}{(\alpha_i,\alpha_i)}= a_{ij}$, $(\delta, \alpha_i) = 0$ for $0\le i,\,j\le n$, and $(\delta, \delta) = 0$. Put $\theta = \delta-\alpha_0$ (cf.~\cite{K}).

Let $\{\,\alpha_i\,|\, i \in I \,\}$ be the set of simple roots for $\mathring{\g}$ and let $\mathring{Q}=\bigoplus_{i \in I}\Z\alpha_i$ be the root lattice of $\mathring{\g}$. Let $\{\,\varpi_i\,|\,i\in I\,\}$ be the fundamental weights for $\mathring{\g}$ and let $\mathring{P}=\bigoplus_{i\in I}\Z\varpi_i$ be the weight lattice of $\mathring{\g}$.
Note that $\theta$ is equal to the maximal root of $\mathring{\g}$ by regarding $\mathring{P}$ as a sublattice of $P / \mathbb{Z}\delta$.

Let $W$ be the affine Weyl group of $\mf{g}$, which is generated by the simple reflection $s_i$ for $0\le i \le n$, and let $\mathring{W}$ be the Weyl group of $\mathring{\mf g}$, which is the subgroup generated by $s_i$ for $i\in I$.
Let $w_0$ be the longest element of $\mathring{W}$.
Note that $W$ is isomorphic to the semidirect product $\mathring{W} \ltimes \mathring{Q}$ under the identification $s_0 \leftrightarrow (s_{\theta}, \theta)$ and $s_i \leftrightarrow (s_i, 0)$.  %where $\theta$ is the maximal root of $\mathring{\mf{g}}$.

Let $\wtd{W}=\mathring{W} \ltimes \mathring{P}$ be the extended affine Weyl group.
Let $\mathcal{T}$ be the set of bijections $\tau : I\cup\{0\} \rightarrow I\cup\{0\}$ such that $a_{\tau(i)\tau(j)} = a_{ij}$ for all $0\le i, j \le n$. 
It is known that each element $\tau \in \mathcal{T}$ induces a unique automorphism $\psi_{\tau}$ of $W$ such that $\psi_{\tau}(s_i) = s_{\tau(i)}$.
Then $W$ is a normal subgroup of $\wtd{W}$ such that $\mathcal{T} \simeq \wtd{W} \big/ W$
and
$\wtd{W} \simeq \mc{T} \ltimes W$, where the action of $\tau \in \mathcal{T}$ in $W$ is given by $\psi_{\tau}$ (see \cite{Bo} for more details).
An expression for $w \in \wtd{W}$ is called reduced if $w = \tau s_{i_1} \cdots s_{i_{\ell(w)}}$, where $\tau \in \mathcal{T}$ and $\ell(w)$ is minimal. We call such $\ell(w)$ the length of $w \in \wtd{W}$.
For $\lambda \in P$, we denote by $t_{\lambda}$ the element $(1, \lambda) \in \wtd{W}$, which is called the translation by $\lambda$.

Let $q$ be an indeterminate. 
We put 
{\allowdisplaybreaks
\begin{gather*}
[m]_q=\frac{q^m-q^{-m}}{q-q^{-1}}\quad (m\in \Z_+),\quad
[m]_q!=[m]_q[m-1]_q\cdots [1]_q\quad (m\geq 1),\quad [0]_q!=1, \\
\begin{bmatrix} m \\ k \end{bmatrix}_q = \frac{[m]_q[m-1]_q\cdots [m-k+1]_q}{[k]_q}\quad (0\leq k\leq m).
\end{gather*}}
If there is no confusion, then we often write $[m]$ and $\begin{bmatrix} m \\ k \end{bmatrix}$ instead of  $[m]_q$ and $\begin{bmatrix} m \\ k \end{bmatrix}_q$ for simplicity, respectively.

\subsection{Quantum loop algebra}
Let ${\bf k}=\C(q)$ be the base field.
The {\em quantum loop algebra} $U_q(\mf{g})$ is the associative $\bf{k}$-algebra generated by $e_i$, $f_i$, and $k^{\pm 1}_i$ for $0\le i\le n$, and $C^{\pm \frac{1}{2}}$ subject to the following relations:
{\allowdisplaybreaks
\begin{gather*}
\text{$C^{\pm \frac{1}{2}}$ are central with $C^{\frac{1}{2}}C^{-\frac{1}{2}}=C^{-\frac{1}{2}}C^{\frac{1}{2}} = 1$,} \\
k_ik_j=k_jk_i,\quad k_ik_i^{-1}=k_i^{-1}k_i=1,\quad \prod_{i=0}^n k_i^{\pm a_i} = ( C^{\pm\frac{1}{2}} )^2, \\
k_ie_jk_i^{-1}=q^{a_{ij}}e_j,\quad k_if_jk_i^{-1}=q^{-a_{ij}}f_j,\\
e_if_j-f_je_i=\delta_{ij}\frac{k_i-k_i^{-1}}{q-q^{-1}},\\
\sum_{m=0}^{1-a_{ij}}(-1)^{m}e_i^{(1-a_{ij}-m)}e_je_i^{(m)}=0,\quad\quad
\sum_{m=0}^{1-a_{ij}}(-1)^{m}f_i^{(1-a_{ij}-m)}f_jf_i^{(m)}=0\quad(i\ne j),
\end{gather*}
for $0\le i,j\le n$, where $e_i^{(m)}=e_i^m/[m]!$ and $f_i^{(m)}=f_i^m/[m]!$ for $0\le i\le n$ and $m\in \Z_+$. 
Here $a_0, a_1, \dots, a_n$ are the numerical labels of the Dynkin diagram associated with $A$ in \cite{K}.
There is a Hopf algebra structure on $U_q(\mf{g})$, where the comultiplication $\Delta$ and the antipode $S$ are given by 
\begin{gather*} 
\Delta(k_i)=k_i\otimes k_i, \quad 
\Delta(e_i)= e_i\ot 1 + k_i\ot e_i,\quad 
\Delta(f_i)= f_i\ot k_i^{-1} + 1\ot f_i,\\  
S(k_i)=k_i^{-1}, \ \ S(e_i)=-k_i^{-1} e_i, \ \  S(f_i)=-f_i k_i,
\end{gather*} 
for $0\le i\le n$.
It is well-known in \cite{B} that as a Hopf algebra, $U_q(\mf{g})$ is also isomorphic to the algebra generated by 
$x_{i,r}^\pm$ ($i\in I, r\in \Z$), $k_{i}^{\pm 1}$ $(i\in I)$, $h_{i,r}$ ($i\in I, r\in \Z\setminus\{0\}$), and $C^{\pm \frac{1}{2}}$ subject to the following relations:
{\allowdisplaybreaks
\begin{gather*}
\text{$C^{\pm \frac{1}{2}}$ are central with $C^{\frac{1}{2}}C^{-\frac{1}{2}}=C^{-\frac{1}{2}}C^{\frac{1}{2}} = 1$,} \\
k_ik_j=k_jk_i,\quad k_ik_i^{-1}=k_i^{-1}k_i=1,\\ 
k_ih_{j,r}=h_{j,r}k_i,\quad k_ix^{\pm}_{j,r}k_i^{-1}=q^{\pm a_{ij}}x^{\pm}_{j,r},\\
[h_{i,r},h_{j,s}]=\delta_{r,-s}\frac{1}{r}[r a_{ij}]\frac{C^r-C^{-r}}{q-q^{-1}},\\
[h_{i,r},x^{\pm}_{j,s}]=\pm \frac{1}{r}[r a_{ij}] C^{\mp |r|/2}x^{\pm}_{j,r+s},\\
x^\pm_{i,r+1}x^\pm_{j,s}-q^{\pm a_{ij}}x^\pm_{j,s}x^\pm_{i,r+1}
=q^{\pm a_{ij}}x^\pm_{i,r}x^\pm_{j,s+1}-x^\pm_{j,s+1}x^\pm_{i,r},\\
[x^+_{i,r},x^-_{j,s}]=\delta_{i,j}\frac{C^{(r-s)/2}\psi^+_{i,r+s}-C^{-(r-s)/2}\psi^-_{i,r+s}}{q-q^{-1}},\\
\sum_{w\in \mf{S}_m}\sum_{k=0}^m {\small \left[\begin{matrix} m \\ k \end{matrix}\right]}x^\pm_{i,r_{w(1)}}\dots x^\pm_{i,r_{w(k)}}x^\pm_{j,s}x^\pm_{i,r_{w(k+1)}}\dots x^\pm_{i,r_{w(m)}}=0\quad (i\neq j),
\end{gather*}
}

\noindent where $r_1,\dots, r_m$ is any sequence of integers with $m=1-a_{ij}$, $\mf{S}_m$ denotes the group of permutations on $m$ letters, and $\psi^\pm_{i,r}$ is the element determined by the following identity of formal power series in $z$;
\begin{equation}\label{eq:psi generators}
\sum_{r=0}^\infty \psi^\pm_{i,\pm r} z^{\pm r} = k_i^{\pm 1}\exp\left( \pm(q-q^{-1}) \sum_{s=1}^\infty h_{i,\pm s} z^{\pm s} \right).
\end{equation}

\subsection{Category $\mc{O}$}

Let $U_q(\mf{b})$ be the subalgebra of $U_q(\mf{g})$ generated by $e_i$, $k_i^{\pm 1}$ for $0\le i\le n$ and $C^{\pm \frac{1}{2}}$. Let $\mf{t}$ be the subalgebra of $U_q(\mf{b})$ generated by $k_i^{\pm 1}$ for $i\in I$, and let $\mf{t}^\ast = ({\bf k}^\times)^I$ be the set of maps from $I$ to ${\bf k}^\times$ which is a group under pointwise multiplication. 

Let $U_q(\mf{g})^\pm$ (resp. $U_q(\mf{g})^0$) be the subalgebras of $U_q(\mf{g})$ generated by $x_{i,r}^{\pm}$ for $i\in I$ and $r\in \Z$ $\big( \text{resp. $k_i^{\pm 1}$, $\psi^\pm_{i,\pm r}$ for $i\in I$, $r>0$ and $C^{\pm \frac{1}{2}}$} \big)$. Then we have a triangular decomposition
\begin{equation*}
U_q(\mf{g})\cong U_q(\mf{g})^-\ot U_q(\mf{g})^0\ot U_q(\mf{g})^+.
\end{equation*}
If we put $U_q(\mf{b})^+ =U_q(\mf{g})^+\cap U_q(\mf{b})$ and $U_q(\mf{b})^0 =U_q(\mf{g})^0\cap U_q(\mf{b})$, then we have
$U_q(\mf{b})^+=\langle\,x_{i,r}^+\,\rangle_{i\in I, r\ge 0}$ and 
$U_q(\mf{b})^0=\langle\,\psi_{i,r}^+, k_i^{\pm 1}, C^{\pm \frac{1}{2}} \,\rangle_{i\in I, r>0}$.

\begin{rem}
{\em 
Throughout this paper, we assume that $C^{\pm \frac{1}{2}}$ acts trivially on a $U_q({\mf b})$-module.
}
\end{rem}

Let $V$ be a $U_q(\mf{b})$-module. For $\omega \in {\mf t}^*$, we define the weight space of $V$ with weight $\omega$ by
\begin{equation*}
	V_{\omega} = \{\, v \in V \, | \, k_i v = \omega(i)v \  (i \in I) \,\}.
\end{equation*}
We say that $V$ is {\em ${\mf t}$-diagonalizable} if $V = \oplus_{\omega \in {\mf t}^*} V_{\omega}$.

A series ${\bf \Psi} = (\Psi_{i, m})_{i \in I, m \ge 0}$ of elements in ${\bf k}$ such that $\Psi_{i, 0} \neq 0$ for all $i \in I$ is called an {\em $\ell$-weight}. 
We often identify ${\bf \Psi} = (\Psi_{i, m})_{m \ge 0}$ with ${\bf \Psi} = (\Psi_i(z))_{i \in I}$, an $I$-tuple of formal power series, where 
\begin{equation*}
	 \Psi_i(z) = \sum_{m \ge 0} \Psi_{i,m} z^m.
\end{equation*}
We denote by ${\mf t}_{\ell}^*$ the set of $\ell$-weights. Since $\Psi_i(z)$ is invertible, ${\mf t}_{\ell}^*$ is a group under multiplication. Let $\varpi : {\mf t}_{\ell}^* \longrightarrow {\mf t}^*$ be the surjective morphism defined by $\varpi(\Psi)(i) = \Psi_{i,0}$ for $i\in I$.
For ${\bf \Psi} \in {\mf t}_{\ell}^*$, we define the {\em $\ell$-weight space} of $V$ with $\ell$-weight ${\bf \Psi}$ by
\begin{equation*}
	V_{\bf \Psi} = \left\{\, v \in V \, \Big| \, \text{for any $i \in I$ and $m \ge 0$, } \exists\, p_{i,m} \in \mathbb{Z}_{+} \,\, \text{such that} \,\, (\psi_{i, m}^+ - \Psi_{i, m} )^{p_{i,m}} v = 0 \right\}.
\end{equation*}
For ${\bf \Psi} \in {\mf t}_{\ell}^*$, we say that $V$ is of {\em highest $\ell$-weight $\bf \Psi$} if there exists a non-zero vector
 $v \in V$ such that
\vskip 2mm
\begin{center}
	(i) $V = U_q(\mf{b})v$ \,\, (ii) $e_i v = 0$ for all $i \in I$ \,\, (iii) $\psi_{i,m}^+ v = \Psi_{i, m}v$ for $i \in I$ and $m \ge 0$.
\end{center}
\vskip 2mm
A non-zero vector $v \in V$ is called a highest $\ell$-weight vector of the weight $\Psi$ if it satisfies the conditions (ii) and (iii).
There exists a unique irreducible $U_q(\mf{b})$-module of highest $\ell$-weight $\bf \Psi$, which we denote by $L({\bf \Psi})$ \cite[Proposition 3.4]{HJ}.

\begin{df}\cite[Definition 3.7]{HJ} \label{df:prefundamental}
{\em
For $r \in I$ and $a \in \bR^{\times}$,
let $L_{r, a}^{\pm}$ be an irreducible $U_q(\mf{b})$-module of highest weight $\bf\Psi$ given by
\begin{equation*}
	\Psi_i(z) = 
	\left\{
		\begin{array}{ll}
			(1-az)^{\pm1} & \text{if $i=r$}, \\
			1 & \text{if $i \neq r$}.
		\end{array}
	\right.
\end{equation*}
The representations $L_{r, a}^{\pm}$ are called {\em prefundamental representations}.
}
\end{df}

For $i\in I$, let $\ov{\alpha}_i \in {\mf t}^*$ given by $\ov{\alpha}_i(j)= q^{a_{ij}}$ ($j\in I$).
We define a partial order $\le$ on ${\mf t}^*$ by
$\omega' \le \omega$  if and only if $\omega'\omega^{-1}$ is a product of $\ov{\alpha}_i^{-1}$'s.
For $\lambda \in {\mf t}^*$, put $D(\lambda) = \{\, \omega \in {\mf t}^* \, | \, \omega \le \lambda \,\}$.

\begin{df}\cite[Definition 3.8]{HJ}\label{df:category O}
{\em 
The category $\mathcal{O}$ consists of $U_q(\mf{b})$-modules $V$ such that
\begin{enumerate}
	\item[(i)] $V$ is ${\mf t}$-diagonalizable,
	\item[(ii)] $\dim V_{\omega} < \infty$ for all $\omega \in {\mf t}^*$,
	\item[(iii)] there exist $\lambda_1, \dots, \lambda_s \in {\mf t}^*$ such that the weights of $V$ are in $\bigcup_{j=1}^s D(\lambda_j)$.
\end{enumerate} 
}	
\end{df}
The category $\mc{O}$ is a tensor category and the simple objects are given as follows.  

\begin{thm}\cite[Theorem 3.11]{HJ}
For ${\bf \Psi} \in {\mf t}_{\ell}^*$, $L({\bf \Psi})$ is in the category $\mc{O}$ if and only if $\Psi_i(z)$ is rational for all $i\in I$.
\end{thm}
Note that for ${\bf \Psi}, {\bf \Psi}'\in {\mf t}_{\ell}^*$, $L({\bf \Psi}{\bf \Psi}')$ is a subquotient of $L({\bf \Psi})\ot L({\bf \Psi}')$. 
Hence any irreducible $U_q(\mf{b})$-module in $\mc{O}$ is a subquotient of a tensor product of prefundamental representations and $1$-dimensional representations.

\section{Unipotent quantum coordinate ring} \label{sec:UQCR}

\subsection{Canonical basis}
Let $U_q^-(\mf{g})$ be the subalgebra of $U_q(\mf{g})$ generated by $f_i$ for $0\le i\le n$. Recall that $U_q^-(\mf{g})$ has a root space decomposition, that is, $U_q^-(\mf{g})=\bigoplus_{\beta\in -Q_+}U_q(\mf{g})_{\beta}$, where $U_q(\mf{g})_{\beta}=\{\,x\,|\,k_ixk_i^{-1}=q^{(\alpha_i,\beta)}x\ (0\le i\le n)\,\}$. Write ${\rm wt}(x)=\beta$ for $x\in U_q(\mf{g})_{\beta}$.

Let us recall the (dual) canonical basis of $U_q^-(\mf{g})$ (see \cite{Kas91,Kas93} for more details).
Given $0\le i\le n$, there exists a unique ${\bf k}$-linear map $e'_i: U_q^-(\mf{g}) \rightarrow U_q^-(\mf{g})$ such that $e'_i(1)=0$, $e'_i(f_j)=\delta_{ij}$ for $0\le j\le n$, and
\begin{equation} \label{eq:derivation}
e'_i(xy)=e'_i(x)y + q^{({\rm wt}(x),\alpha_i)}xe'_i(y),
\end{equation}
for homogeneous $x,y\in U_q^-(\mf{g})$.
Then there exists a unique non-degenerate symmetric ${\bf k}$-valued bilinear form $(\ ,\ )$ on $U_q^-(\mf{g})$ such that
\begin{equation}\label{eq:bilinear form}
(1,1)=1,\quad (f_ix,y)=(x,e'_i(y)),
\end{equation}
for $0\le i\le n$ and $x,y\in U_q^-(\mf{g})$.

For any homogeneous $x\in U_q^-(\mf{g})$, we have $x=\sum_{k\ge 0}f_i^{(k)}x_k$, where $e'_i(x_k)=0$ for $k\ge 0$. Then we define $\wtd{f}_ix=\sum_{k\ge 0}f_i^{(k+1)}x_k$.
Let $\mc{A}_0$ denote the subring of ${\bf k}$ consisting of rational functions regular at $q=0$.
Let 
\begin{equation*}
\begin{split}
L(\infty)&= \sum_{r\ge 0, 0\le i_1,\dots,i_r\le n}\mc{A}_0\wtd{f}_{i_1}\dots \wtd{f}_{i_r}1,\\
B(\infty)&=\{\,\wtd{f}_{i_1}\dots \wtd{f}_{i_r}1 \!\!\pmod{qL(\infty)} \,|\,r\ge 0, 0\le i_1,\dots,i_r\le n\,\}\setminus\{0\}\subset L(\infty)/qL(\infty).
\end{split}
\end{equation*}
The pair $(L(\infty),B(\infty))$ is called the {\em crystal base} of $U_q^-(\mf{g})$.

Let $\mc{A}=\Z[q,q^{-1}]$, and let $U_q^-(\mf{g})_{\mc{A}}$ be the $\mc{A}$-subalgebra of $U_q^-(\mf{g})$ generated by $f_i^{(k)}$ for $0\le i\le n$ and $k\in \Z_+$. Let $- : U_q^-(\mf{g})\rightarrow U_q^-(\mf{g})$ be the automorphism of $\C$-algebras given by $\ov{q}=q^{-1}$ and $\ov{f_i}=f_i$ for $0\le i\le n$.
Then $(L(\infty),\ov{L(\infty)},U_q^-(\mf{g})_{\mc{A}})$ is balanced, that is, 
the map
\begin{equation}\label{eq:balanced}
\xymatrixcolsep{2pc}\xymatrixrowsep{-0.3pc}\xymatrix{
E:=L(\infty)\cap \ov{L(\infty)}\cap U_q^-(\mf{g})_{\mc{A}} \ \ar@{->}[r]  &\ L(\infty)/qL(\infty) \\
	x \ \ar@{|->}[r] & \, x \,\, \scalebox{0.8}{\text{(mod $qL(\infty)$)}}
}
\end{equation}
is a $\C$-linear isomorphism.

Let $G$ denote the inverse of the map \eqref{eq:balanced}. 
Then 
\begin{equation*}
G(\infty):=\{\,G(b)\,|\,b\in B(\infty)\,\}
\end{equation*} 
is an $\mc{A}$-basis of $U_q^-(\mf{g})_{\mc A}$, which is called the {\em canonical basis} or {\em global crystal basis}.
Let
\begin{equation*}
G^{\rm up}(\infty)=\{\,G^{\rm up}(b)\,|\,b\in B(\infty)\,\}
\end{equation*}
be the dual basis of $G(\infty)$ with respect to the bilinear form \eqref{eq:bilinear form}, that is, 
$(G^{\rm up}(b),G(b'))=\delta_{bb'}$ for $b,b'\in B(\infty)$.
We call $G^{\rm up}(\infty)$ the {\em dual canonical basis} of 
$U_q^-(\mf{g})^{\rm up}_{\mc{A}}:=\{\,x\in U_q^-(\mf{g})\,|\,(x,U_q^-(\mf{g})_{\mc{A}})\in \mc{A}\,\}$.

Let $U_q(\mf{g})^e$ be the ${\bf k}$-algebra generated by $U_q(\mf{g})$ and $q^{\pm d}$ where $q^d$ commutes with $k_i^{\pm 1}$ and satisfies $q^{d}q^{-d}=1$, $q^de_iq^{-d}=q^{\delta_{0 i}}e_i$, and $q^df_iq^{-d}=q^{-\delta_{0 i}}e_i$.
Given a dominant integral weight $\La$ for $\mf{g}$, let $V(\La)$ be the irreducible highest weight module over $U_q(\mf{g})^e$. Let $B(\La)$ and $G(\La)=\{\,G_\La(b)\,|\,b\in B(\La)\,\}$ denote the crystal and canonical basis of $V(\La)$, respectively. Let $\ast$ be the ${\bf k}$-algebra anti-automorphism of $U_q^-(\mf{g})$ given by $(f_i)^\ast=f_i$ for $0\le i\le n$. Then $L(\infty)^\ast=L(\infty)$ and $B(\infty)^\ast=B(\infty)$.
We may regard $B(\La)\subset B(\infty)$ by
\begin{equation}\label{eq:B(Lambda)}
B(\La)=\{\,b\in B(\infty)\,|\,\varepsilon_i^\ast(b)\leq \langle \La,\alpha_i^\vee\rangle \ (0\le i\le n)\,\},
\end{equation}
where $\alpha_i^\vee$ is the simple coroot and $\varepsilon_i^\ast(b)=\max\{\,k\,|\,b^\ast=\wtd{f}_i^k (b_0^\ast)\ \text{for some $b_0\in B(\infty)$}\,\}$.
We have $G_\La(b)=\pi_\La(G(b))$ for $b\in B(\Lambda)$, where $\pi_\La : U_q^-(\mf{g})\longrightarrow V(\La)$ is the canonical projection.

Let $G^{\rm up}(\La)=\{\,G^{\rm up}(b)\,|\,b\in B(\La)\,\}$ be the dual basis of $G(\La)$ with respect to the bilinear form on $V(\La)$ in \cite[(4.2.4), (4.2.5)]{Kas93}. 
Let 
\begin{equation*}
\xymatrixcolsep{2pc}\xymatrixrowsep{4pc}\xymatrix{
\iota_\La : V(\La)^\vee \ar@{->}[r]  & U_q(\mf{g})^\vee.}
\end{equation*}
be the dual of $\pi_\La$, where $V(\La)^\vee$ and $U_q^-(\mf{g})^\vee$ are the duals of $V(\La)$ and $U_q^-(\mf{g})$ with respect to \cite[(4.2.4), (4.2.5)]{Kas93} and \eqref{eq:bilinear form}, respectively.
Then we have 
\begin{equation}\label{eq:embedding and dual canonical basis}
\iota_\La(G^{\rm up}_\La(b))=G^{\rm up}(b),
\end{equation}
for $b\in B(\La)$. Here we identify $V(\La)^\vee$ with $V(\La)$, and $U_q^-(\mf{g})^\vee$ with $U_q^-(\mf{g})$.

\subsection{Unipotent quantum coordinate ring $U_q^-(w)$}
Let us recall the unipotent quantum coordinate ring (cf.~\cite{Ki12})

For $0\le i\le n$, let $T_i$ be the ${\bf k}$-algebra automorphism of $U_q(\mf{g})$ \cite{Lu93} given by
\begin{gather*}
T_i(k_j) = k_jk_i^{-a_{ji}}, \\ 
T_i(e_i) = -f_i k_i,\quad
T_i(e_j) = \sum_{r+s=-a_{ij}} (-1)^r q^{-r} e_i^{(s)} e_j e_i^{(r)}\quad (j\neq i), \\
T_i(f_i) = -k_i^{-1}e_i, \quad 
T_i(f_j) = \sum_{r+s=-a_{ij}} (-1)^r q^{ r} f_i^{(r)} f_j f_i^{(s)}\quad (j\neq i),
\end{gather*}
for $0\le j\le n$. Note that $T_i=T''_{i,1}$ in \cite{Lu93}.

Let $\Delta^{\pm}$ be the set of positive (resp. negative) roots of $\mf{g}$. For $w\in W$, let $\Delta^+(w)=\Delta^+\cap w\Delta^-$. 
Suppose that the length of $w$ is $\ell$.
Let $R(w)=\{\,\wtd{w}=(i_1,\dots,i_\ell) \,\,|\,\, 0 \le i_j \le n \,\,\, \text{and} \,\, w=s_{i_1}\dots s_{i_\ell}\,\}$ be the set of reduced expressions of $w$. 
A $2$-braid move on $\wtd{w} \in R(w)$ is defined by $(\dots, i, j, \dots) = (\dots, j, i, \dots)$ for $i,\, j \in I$ such that $|i-j| > 1$.
For $\wtd{w}=(i_1,\dots,i_\ell)\in R(w)$, we have $\Delta^+(w)=\{\,\beta_k\,|\,1\le k\le \ell\,\}$, where 
\begin{equation}\label{eq:root vector}
\beta_k=s_{i_1}\dots s_{i_{k-1}}(\alpha_k)\quad (1\le k\le \ell).
\end{equation}

For $1\le k\le \ell$ and $c\in \Z_+$, let 
$F(c\beta_k)=T_{i_1}\dots T_{i_{k-1}}(f_{i_k}^{(c)})\in U_q^-(\mf{g})_{-c\beta_k}$. 
For ${\bf c}=(c_1,\dots,c_\ell)\in \Z_+^\ell$, put
\begin{equation*}
F({\bf c},\wtd{w})=F(c_1\beta_1) \cdots F(c_\ell\beta_\ell).
\end{equation*}
Then $\{\,F({\bf c},\wtd{w})\,|\,{\bf c}\in \Z_+^\ell\,\}\subset U_q^-(\mf{g})_{\mc A}$. 
Furthermore, 
we have $F({\bf c},\wtd{w})\in L(\infty)$ and $F({\bf c},\wtd{w}) \!\!\pmod{qL(\infty)}$ belongs to $B(\infty)$, say $b({\bf c},\wtd{w})$ \cite{S94}.

We assume that $\Delta^+(w)$ is linearly ordered by $\beta_1<\dots<\beta_\ell$.
Then the following $q$-commutation relation for the root vectors $F(c_k\beta_k)$'s holds \cite{LS91}:
\begin{equation}\label{eq:LS formula}
F(c_j\beta_j)F(c_i\beta_i)-q^{-(c_i\beta_i,c_j\beta_j)}F(c_i\beta_i)F(c_j\beta_j)
=\sum_{{\bf c}'}f_{{\bf c}'}F({\bf c}',\wtd{w}),
\end{equation}
for $i<j$ and $c_i,c_j \in \Z_+$, where the sum is over ${\bf c}'=(c'_k)$ such that $c_i\beta_i+c_j\beta_j=\sum_{i\le k \le j}c'_{k}\beta_k$ with $c'_i<c_i$, $c'_j<c_j$ and $f_{{\bf c}'} \in {\bf k}$.

\begin{df}
{\em 
For $w \in W$, we denote by $U_q^-(w)$ the vector space over ${\bf k}$ generated by $\{\,F({\bf c},\wtd{w})\,|\,{\bf c}\in \Z_+^\ell\,\}$.
}
\end{df}
Note that $U_q^-(w)$ does not depend on the choice of $\wtd{w}\in R(w)$, and the $q$-commutation relation \eqref{eq:LS formula} implies that $U_q^-(w)$ is the ${\bf k}$-subalgebra of $U_q^-(\mf{g})$ generated by $\{\,F(\beta_k)\,|\,1\le k\le \ell\,\}$.
Moreover, it is shown in \cite[Theorem 4.25]{Ki12} that $U_q^-(w)$ is compatible with the dual canonical basis in the following sense
\begin{equation}\label{eq:dual canonical and quantum nilpotent}
U_q^-(w)\cap U_q^-(\mf{g})_{\mc{A}}^{\rm up} = \bigoplus_{b\in B(w)}\mc{A}G^{\rm up}(b),
\end{equation} 
where $B(w)=\{\,b({\bf c},\wtd{w})\,|\,{\bf c}\in \Z_+^\ell\,\}$.
Let $\{\,F^{\rm up}({\bf c},\wtd{w})\,|\,{\bf c}\in \Z_+^\ell\,\}$ be the dual basis of $\{\,F({\bf c},\wtd{w})\,|\,{\bf c}\in \Z_+^\ell\,\}$. Then the left-hand side of \eqref{eq:dual canonical and quantum nilpotent} is also the $\mc{A}$-subalgebra generated by $\{\,F^{\rm up}(\beta_k)\,|\,1\le k\le \ell\,\}$.
Note that the same formula \eqref{eq:LS formula} also holds for the dual root vectors 
$F^{\rm up}(\beta_k)$ where $f_{{\bf c}'}\in \mc{A}$ (see \cite[Theorem 4.27]{Ki12}).

The subalgebra $U_q^-(w)$ is called the {\em unipotent quantum coordinate ring} since the quotient $\C\ot_{\mc{A}}U_q^-(w)$ with respect to the $\mc{A}$-lattice generated by $F^{\rm up}({\bf c},\wtd{w})$ is isomorphic to the coordinate ring of the unipotent subgroup $N(w)$ of the Kac-Moody group associated to $w$ \cite[Theorem 4.44]{Ki12}.

\subsection{$U_q(\mf{b})$-action on $U_q^-(\varpi_r)$} \label{subsec:b-actions}
From now on, we assume that $\mf g$ is of type $A_n^{(1)}$ $(n \ge 1)$ or $D_n^{(1)}$ $(n \ge 4)$ and $r \in I$ denotes an index such that $\varpi_r$ is minuscule or the fundamental representation of $\mathring{\mf g}$ with highest weight $\varpi_r$ is minuscule, that is,
\begin{equation}\label{eq:minuscule index}
		\begin{cases} 
			r \in \{\,1,\dots,n\,\} & \text{for type $A_n$}, \\
			r \in\{\,1,n-1,n\,\} & \text{for type $D_n$}.
		\end{cases}
\end{equation}

\begin{lem}\label{lem:property of w_r}
Let $t_{\varpi_r}$ be the translation by $\varpi_r$.
Suppose that
$$t_{\varpi_r}^{-1}  = w_r \tau,$$ 
for some $\tau\in\mc{T}$ and $w_r\in W$ of length $\ell$.
Then $w_r$ satisfies the following properties:
\vskip 2mm

\begin{enumerate}
\item $w_r$ is the maximal element in the set $\mathring{W}_{r}\setminus\mathring{W}$ of minimal length coset representatives, where $\mathring{W_r}$ is the subgroup of stabilizers of $\varpi_r$ in $\mathring{W}$,
\vskip 1mm

\item a reduced expression of $w_r$ is unique up to $2$-braid moves,
\vskip 1mm

\item $\Delta^+(w_r) = \{\,  \beta \in \Delta^+ \cap \mathring{Q} \mid (\varpi_r, \beta) = 1 \,\}$,
\vskip 1mm

\item if $\beta_k$ $(1\le k\le \ell)$ is a positive root in \eqref{eq:root vector} with respect to a reduced expression of $w_r$, then we have $\beta_1 = \alpha_r$ and $\beta_\ell=\theta$,
\vskip 1mm

\item $w_r^{-1} = w_{r^*}$, where $r^*$ is determined by $\alpha_{r^*} = -w_0(\alpha_r)$.
\end{enumerate}
\end{lem}
\pf 
Let  $J = I \setminus \{ r \}$ (resp.\ $J^* = I \setminus \{ r^* \}$) and $w_0^J$ (resp. $w_0^{J^*}$) be the longest element of the parabolic subgroup $\mathring{W}_J$ (resp.\ $\mathring{W}_{J^*}$) of $\mathring{W}$.
The assertion (1) is well-known (cf.~\cite{Bo}). Thus we have 
\begin{align} \label{Eq: w_0 J}
w_0 = w_0^J w_r = w_r w_0^{J^*}.
\end{align}

Since $w_r^{-1}$ is $\varpi_r$-minuscule (see \cite[(2.1)]{Stem} for its definition),
 (2) follows from \cite[Proposition 2.1]{Stem} (see also \cite[Remark 5.2]{JK20}).
 By the definition of being $\varpi_r$-minuscule, 
 we have 
 $$
 \Delta^+(w_r) \subset \{\,  \beta \in \Delta^+ \cap \mathring{Q} \mid (\varpi_r, \beta)=1 \,\}.
 $$
 As $ |\Delta^+(w_r)| =  \ell(w_r) =  \ell( w_0^J w_0  ) = \ell(w_0) - \ell(w_0^J) $, one can check that 
 $$
  | \Delta^+(w_r) | = \# \{\,  \beta \in \Delta^+ \cap \mathring{Q} \mid (\varpi_r, \beta)=1 \,\},
  $$
  which gives (3).

We take a reduced expression $\wtd{w}_r = s_{i_1} s_{i_2}  \cdots s_{i_\ell} $ of $w_r$. 
It is easy to see that $ w_0^J (\alpha_r) = \theta $.
It follows from (1) and \cite[Lemma 2.2.1]{GP} that $i_1  =  r$ and $i_\ell = r^*$.
 Thus $\beta_1 = \alpha_r$ and 
$$
\beta_\ell = - w_r (\alpha_{r^*})  =  w_0^J w_0 ( - \alpha_{r^*}) = w_0^J ( \alpha_r) = \theta.
$$
Since $w_0$ and $w_0^{J^*}$ are involutions, it follows from \eqref{Eq: w_0 J} that 
$$ 
w_{r^*} = w_0^{J^*} w_0  = w_r^{-1},
$$
which implies (5).
\qed

\begin{rem}\label{rem:formula for w_i}
{\rm
For $i \leq j$, put
$$
	{\bf s}_{(i,\,j)} = s_i s_{i+1} \dots s_j.
$$

A reduced expression of $w_r$ is given as follows:
\begin{equation*}
	w_r =
	\begin{cases}
		{\bf s}_{(r,\,n)} {\bf s}_{( r-1,\,n-1)} \dots {\bf s}_{( 1,\,n-r+1)} & \text{for type $A_n^{(1)}$ with $r \in I$,} \\
		{\bf s}_{(1,\,n)} {\bf s}_{(1,\,n-2)}^{-1} & \text{for type $D_n^{(1)}$ with $r = 1$,}\\
		\tilde{\bf s}_1 \tilde{\bf s}_2 \dots \tilde{\bf s}_{n-1} & \text{for type $D_n^{(1)}$ with $r = n$,}
	\end{cases} 
\end{equation*}
where $\tilde{\bf s}_k$ is given by
\begin{equation*}
	\tilde{\bf s}_k =
	\begin{cases}
		s_n {\bf s}_{(k,\, n-2)}^{-1} & \text{if $k$ is odd,} \\
		{\bf s}_{(k,\, n-1)}^{-1} & \text{if $k$ is even,} \\
		s_n & \text{if $n$ is even and $k = n-1$.}
	\end{cases}
\end{equation*}
Note that for type $D_n^{(1)}$, a reduced expression of $w_{n-1}$ is obtained from the one of $w_n$ by replacing $s_{n}$ with $s_{n-1}$.
}
\end{rem}

By abuse of notation, we put 
\begin{equation} \label{eq:definition of Uwr}
U_q^-(\varpi_r):= U_q^-(w_r).
\end{equation}
The following is the main result in this section.
\begin{thm}\label{thm:main-1}
For a minuscule $\varpi_r$, there is a $U_q(\mf{b})$-module structure on $U_q^-(\varpi_r)$ such that $U_q^-(\varpi_r)$ is in the category $\mc{O}$.
\end{thm}

The rest of this subsection is devoted to proving Theorem \ref{thm:main-1}. 
To define a $U_q(\mf{b})$-action on $U_q^-(\varpi_r)$, we will consider the ${\bf k}$-linear operators \eqref{eq:borel action} on $U_q^-(\mf{g})$, and then verify that the operators satisfy the defining relations of $U_q(\mf{b})$.
\vskip 2mm

Let us choose $\wtd{w}=(i_1,\dots,i_\ell)\in R(w_r)$ and let $\{\,\beta_k\,|\,1\le k\le \ell\,\}$ be as in \eqref{eq:root vector} with respect to $\wtd{w}$. 
By Lemma \ref{lem:property of w_r}(2), $F(\beta_k)$ (and hence $F^{\rm up}(\beta_k)$) is independent of the choice of $\wtd{w}$ up to permutations. For $1\le k\le \ell$, put ${\bf 1}_{\beta_k} \in B(\infty)$ such that
\begin{equation*}
{\bf 1}_{\beta_k} \equiv F^{\rm up}(\beta_k) \!\!\pmod{qL(\infty)} \in B(\infty).
\end{equation*}
Let ${\rm wt} : B(\infty) \, \rightarrow P$ denote the weight function. For example, ${\rm wt}({\bf 1}_{\beta_k}) = -\beta_k$.

\begin{rem}\label{rem:minuscule case}
{\rm
Note that for $b \in B(\La)$ and $m \in \Z_{>0}$, we have 
\begin{equation}\label{eq:useful fact 1}
\begin{split}
& \te_i^{\,m}(b) = {\bf 0} \,\, \Longleftrightarrow \,\,
e_i^m G^{\rm up}_\La(b) = 0 \,\, \Longleftrightarrow \,\,
(e_i')^m G^{\rm up}(b) = 0
\end{split}
\end{equation}
(\cite[Lemma 5.1.1]{Kas93}).
Also, for $b\in B(\La)$ with $\vep_i(b)=1$, we have
\begin{equation}\label{eq:useful fact 2}
G^{\rm up}_\La(\te_i(b)) = e_iG^{\rm up}_\La(b), 
\end{equation}
which is equivalent to $G^{\rm up}(\te_i(b)) = e_i'G^{\rm up}(b)$ (cf. \cite[Theorem 3.14]{Ki12}).}
\end{rem}

\begin{lem} \label{lem:dual PBW properties 1}
Let $1 \le k \le \ell$ be given.

\begin{enumerate}
	\item We have 
$$F^{\rm up}(\beta_k)\in \iota_{\La_r}(G^{\rm up}(\La_r)),$$
where $\La_r$ is the $r$-th fundamental weight for $\mf{g}$. 
In particular, ${\bf 1}_{\beta_k}\in B(\La_r)$ by regarding $B(\Lambda_r)$ as a subset of $B(\infty)$ $($cf.~\eqref{eq:B(Lambda)}, \eqref{eq:embedding and dual canonical basis}$)$.
	\vskip 1mm
	
	\item For $i\in I$, we have 
$$
 e'_i ( F^{\rm up}(\beta_k) ) = 
 \begin{cases}
 	 F^{\rm up}(\beta_k - \alpha_i)  & \text{ if }  (\alpha_i, \beta_k) = 1 + \delta_{i,r},  \\
 	 0 & \text{ otherwise.}
 \end{cases}
$$
Here we understand $F^{\rm up}(0) = 1$.
\end{enumerate}
\end{lem}
\pf 
%%%%%%
(1) First, we have $F^{\rm up}(\beta_k)\in G^{\rm up}(\infty)$ for $1\le k\le \ell$ by \cite[Theorem 4.29]{Ki12}. So by \eqref{eq:embedding and dual canonical basis}, it remains to show that ${\bf 1}_{\beta_k}\in B(\La_r)$, equivalently 
\begin{equation}\label{eq:varepsilon*}
\vep_i^*({\bf 1}_{\beta_k}) = \delta_{i,r}.
\end{equation}

It is clear that $\vep_0^*({\bf 1}_{\beta_k})=0$ since $\beta_k-\alpha_0\not\in Q_+$.
Let $i \in I\setminus\{r\}$ be given.
By Lemma \ref{lem:property of w_r}(1), we have 
		\begin{equation*}
			\begin{split}
			\ell(s_i s_{i_1}\dots s_{i_{k-1}}) &= \ell(s_{i_1}\dots s_{i_{k-1}}) + 1, \\
			\ell(s_i s_{i_1}\dots s_{i_{k}}) &= \ell(s_{i_1}\dots s_{i_{k}}) + 1,
			\end{split}
		\end{equation*}
for $1 \le k \le \ell$. 
Then we have $T_i F^{\rm up}(\beta_k) \in U_q^-(\mf{g})$ \cite{Lu93} (cf. \cite[Propositions 7.1 and 7.4]{GLS13}), and hence $\varepsilon_i^* ({\bf 1}_{\beta_k}) = 0$ since $T_i^{-1}(U_q^-(\mf{g}))\cap U_q^-(\mf{g})={\rm Ker}(e'_i)^\ast$ (see \cite[(3.4.4)]{S94}).
On the other hand, we have  $\varepsilon_r^* ({\bf 1}_{\beta_k}) \neq 0$ since ${\rm wt}( {\bf 1}_{\beta_k}) \ne 0$. Thus we conclude from the weight consideration that $\varepsilon_r^* ({\bf 1}_{\beta_k}) = 1$.
This proves that  ${\bf 1}_{\beta_k}\in B(\La_r)$.
%%%%%%
\vskip 2mm
%%%%%%
(2) 
Let $B_0(\La_r)$ be the connected component of $1$ in $B(\La_r)$ generated by $\wtd{e}_i, \wtd{f}_i$ for $i\in I$. We denote by $\mathring{W} \cdot \varpi_r$ the orbit of $\varpi_r$ by the action of $\mathring{W}$.
Note that $B_0(\La_r)$ is the crystal of the fundamental representation $V(\varpi_r)$ of $U_q(\mathring{\mf g})$ with highest weight $\varpi_r$. Thus, the map 
\begin{equation*}
\xymatrixcolsep{2pc}\xymatrixrowsep{-0.3pc}\xymatrix{
B_0(\La_r) \ \ar@{->}[r]  &\ \mathring{W} \cdot \varpi_r \\
	b \ \ar@{|->}[r] & \, {\rm wt}(b) 
}
\end{equation*}
is bijective because $B_0(\La_r)$ is the crystal of the minuscule representation $V(\varpi_r)$.
By (1), we have 
\begin{equation} \label{eq:weights in B0}
 \{\,  \varpi_r - \beta_k \mid k=1,\ldots, \ell \, \} \subset  \mathring{W} \cdot \varpi_r.
\end{equation}
By \eqref{eq:weights in B0} and \cite[Lemma 1.4]{AK}, one can check that 
\begin{equation*}
\varpi_r - \beta_k + \alpha_i \in \mathring{W} \cdot \varpi_r 
\,\,\,
\Longleftrightarrow
\,\,\,	
( \varpi_r - \beta_k,  \alpha_i) = -1
\,\,\,
\Longleftrightarrow
\,\,\,
(\alpha_i, \beta_k) = 1 + \delta_{i,r}.
\end{equation*}
 This yields that 
\begin{equation} \label{eq:epsilon i on beta}
\varepsilon_i ( {\bf 1}_{\beta_k} ) = 
\begin{cases}
	1  & \text{ if }(\alpha_i, \beta_k) = 1 + \delta_{i,r},  \\
	0 & \text{ otherwise.}
\end{cases}
\end{equation}
Therefore we have the assertion by  \eqref{eq:useful fact 1} and \eqref{eq:useful fact 2}.
\qed

\begin{rem}
{\rm
One can also check \eqref{eq:varepsilon*} directly by the formula in \cite[Theorem 3.7]{BZ} or \cite[Proposition 7.4]{Re}. For example, see \cite[Theorem 3.12]{JK19} in which the authors give a combinatorial description of $\varepsilon_n^*$ for type $D_n$ by using \cite[Theorem 3.7]{BZ}.
}
\end{rem}

Consider the following root vector with the weight $\theta$;
\begin{equation} \label{eq:xz}
\mathsf{x}_0:=F^{\rm up}(\beta_\ell)=F^{\rm up}(\theta).
\end{equation}
Let $J_0= \left\{ \, i\in I \, | \, (\al_i, \theta)=0 \, \right\}$ and $J_1= I \, \setminus \, J_0$.
If ${\mf g}$ is of type $A_n^{(1)}$, then we have
\begin{equation*}
	J_1 =
	\begin{cases}
		\{ \, 1 \, \} & \text{if $n \ge 1$,} \\
		\{ \, 1,\, n \, \} & \text{if $n \ge 2$,}
	\end{cases}
	\qquad\,\,
	(\theta, \alpha_i) = 
	\begin{cases}
		2 & \text{if $n = 1$ and $i \in J_1$,} \\
		1 & \text{if $n \ge 2$ and $i \in J_1$.}
	\end{cases}
\end{equation*}
If ${\mf g}$ is of type $D_n^{(1)}$, then we have $J_1 = \{ \, 2 \, \}$ and $(\theta, \alpha_2) = 1$.

\begin{lem}\label{lem:dual PBW properties 2} \
\begin{itemize}
\item[(1)] For $i\in I$, we have
	\begin{equation*}
		\left\{
			\begin{array}{ll}
				e_i'(\xz) = 0 & \text{if $i \in J_0$}, \\
				(e_i')^2 (\xz) = 0 & \text{if $i \in J_1$}.
			\end{array}
		\right.		
	\end{equation*}
	
\item[(2)] For $i\in J_1$, we have
\begin{equation*}
	\xz \cdot e_i'(\xz) = q^{-(1-\delta_{1,n})} e_i'(\xz) \cdot \xz\,.
\end{equation*}
\end{itemize}
\end{lem}
\pf 
Note that (1) and (2) are straightforward for type $A_1^{(1)}$, so we assume that $n \ge 2$ when ${\mf g}$ is of type $A_n^{(1)}$.
\vskip 2mm

(1) By Lemma \ref{lem:dual PBW properties 1}, ${\bf 1}_{\theta}={\bf 1}_{\beta_\ell}$ belongs to $B(\La_r)$. 
Let us consider two cases.

\vskip 3mm

{\em Case 1}. $r \in I \, \setminus \, \{ 1,\, n\}$ in type $A_n^{(1)}$ or $r \in \{ 1,\, n-1,\, n \}$ in type $D_n^{(1)}$.
In this case, 
by \eqref{eq:epsilon i on beta}, one can check that
\begin{equation} \label{eq:value of bilinear form}
(\alpha_i,\, \theta) = 1 + \delta_{i, r} \, \Longleftrightarrow \, i \in J_1.
\end{equation}
Thus, we have
\begin{equation*}\label{eq:varespilon beta_k}
	\varepsilon_i ({\bf 1}_{\theta}) =
		\begin{cases} 
			0 & \text{if $i \in J_0$}, \\
			1 & \text{if $i \in J_1$}.
		\end{cases} 
\end{equation*}

For $i\in J_0$, we have $\wtd{e}_i{\bf 1}_{\theta}=0$ in $B(\La_r)$, which is equivalent to $\wtd{e}_i {\bf 1}_{\theta} =0$ in $B(\infty)$. Then it is equivalent to $e'_i(\xz)=0$ (see \eqref{eq:useful fact 1}).
For $i\in J_1$, we have $\varepsilon_i({\bf 1}_\theta)=1$, and 
\begin{align*} \label{eq:commuting when varepsilon is 1}
&e_i'\xz=e_i'G^{\rm up}({\bf 1}_\theta) = G^{\rm up}( \wtd{e}_i{\bf 1}_\theta),\\
&(e_i')^2\xz=(e_i')^2\,G^{\rm up}({\bf 1}_\theta) = G^{\rm up}\left( \widetilde{e}^2_i{\bf 1}_\theta\right)=0,
\end{align*}
(see \eqref{eq:useful fact 1} and \eqref{eq:useful fact 2}).

\vskip 2mm

{\em Case 2}. $r = 1$ or $n$ in type $A_n^{(1)}$.
We may assume that $r = 1$, since the proof of the case of $r = n$ is almost identical.
In this case, since
\begin{equation} \label{eq:value of bilinear form 2}
(\alpha_i,\, \theta) = 1 + \delta_{i, 1} \, \Longleftrightarrow \, i = n,
\end{equation}
we have
\begin{equation*}
	\varepsilon_i ({\bf 1}_{\theta}) =
		\begin{cases} 
			0 & \text{if $i \neq n$}, \\
			1 & \text{if $i = n$}.
		\end{cases} 
\end{equation*}
By the same argument as in {\em Case 1}, we conclude that $e_i' \xz = 0$ for $i \neq n$, and $(e_n')^2 \xz = 0$.
Note that it is obvious that $(e_1')^2 \xz = 0$, that is, $(e_i')^2 \xz = 0$ for $i \in J_1$.
\vskip 1mm

Hence, we complete the proof of (1).
\vskip 2mm

(2)
Put 
\begin{equation*}
	\tilde{J}_1 = 
	\begin{cases}
		J_1 \, \setminus \, \{ \delta_{1,r}\cdot 1,\, \delta_{n,r}\cdot n \} & \text{if ${\mf g}$ is of type $A_n^{(1)}$}, \\
		J_1 & \text{if ${\mf g}$ is of type $D_n^{(1)}$.}
	\end{cases}
\end{equation*}
By Lemma \ref{lem:dual PBW properties 1}(2) and \eqref{eq:value of bilinear form}-\eqref{eq:value of bilinear form 2}, we have
$\wtd{e}_i{\bf 1}_\theta={\bf 1}_{\theta - \alpha_i}$ for $i \in \tilde{J}_1$.
Then one can check the required relation directly by using the formula \eqref{eq:LS formula}. We remark that for the case of type $A_n^{(1)}$ with $r = 1$ or $n$, since $e_{ \delta_{1,r}\cdot 1 + \delta_{n,r}\cdot n}' (\xz)  = 0$, we also obtain the required relation in this case.
\qed

\vskip 2mm

Let $a \in \bR^\times$.
For $0\le i\le n$, let $\pfe_i$ and $\pfk_i$ be the ${\bf k}$-linear operators on $U_q^-(\mf{g})$ given by
\begin{equation} \label{eq:borel action}
\mathsf{e}_i(u)= 
\begin{cases}
e_i' (u) \quad &\text{ if } i\in I\,,  \\
a q^{-(\theta,\, \beta)} u \cdot \xz \quad &\text{ if } i=0\,,
\end{cases}
\qquad
\mathsf{k}_i(u)=  
\begin{cases}
q^{(\al_i,\, \beta)}u \quad &\text{ if } i\in I\,,  \\
q^{ - (\theta,\, \beta)}u  \quad &\text{ if } i=0\,.
\end{cases} 
\end{equation}
for $u\in U_q^-(\mf{g})_\beta$ ($\beta\in -Q_+$).

\begin{ex} \label{ex:A1(1)}
{\em 
Let us consider the above ${\bf k}$-linear operators on $U_q^-({\mf g})$ for type $A_1^{(1)}$.
In this case, we get
\begin{equation*}
	w = s_1, \quad \theta = \alpha_1, \quad \xz = f_1.
\end{equation*}

By \eqref{eq:derivation}, one can check that
\begin{equation} \label{eq:action of e' on fm}
	e_1' \left( f_1^m \right) = q^{-m+1} [m] f_1^{m-1},
\end{equation}
where we understand $f_1^{-1} = 0$ and $f_1^0 = 1$ (cf. \cite[Section 3.1]{Kas91}).
\vskip 2mm

Put 
$S_{i, j} = \pfe_i^3 \pfe_j - [3] \pfe_i^2 \pfe_j \pfe_i + [3] \pfe_i \pfe_j \pfe_i^2 - \pfe_j \pfe_i^3$.
Here $[3] = q^{-2} + 1 + q^2$ by definition.
We show that 
\begin{equation*}
	S_{i,j}(u) = 0,
\end{equation*}
for $u \in U_q^-({\mf g})_\beta$. Put $s = (\theta, \beta)$ and $t = (\beta, \alpha_1)$.
\vskip 2mm

{\it Case 1}. $(i, j) = (0, 1)$. 
By definition \eqref{eq:borel action}, one can check that
\begin{equation*}
\begin{split}
	\pfe_0^3 \pfe_1 (u) &= a^3 q^{-3s} e_1'(u) f_1^3\,, \\
	\pfe_0^2 \pfe_1 \pfe_0 (u) &= a^3 q^{-3s+2} e_1'(u) f_1^3 + a^3 q^{-3s+t+2} u f_1^2\,, \\
	\pfe_0 \pfe_1 \pfe_0^2 (u) &= a^3 q^{-3s+4} e_1'(u) f_1^3 + a^3 q^{-3s+t+2} (1+q^2) u f_1^2\,, \\
	\pfe_1 \pfe_0^3 (u) &= a^3 q^{-3s+6} e_1'(u) f_1^3 + a^3 q^{-3s+t+2} (1 + q^2 +  q^4 ) u f_1^2\,.
\end{split}
\end{equation*}
Thus, we have $S_{0,1}(u) = 0$ in this case.

\vskip 2mm

{\it Case 2}. $(i,j) = (1,0)$.  Similarly, we have
\begin{equation*}
\begin{split}
	\pfe_0^3 \pfe_1 (u) &= aq^{-s} e_1'^3(u)f_1 + aq^{-s+t} (1 + q^2 + q^4) e_1'^2(u)\,, \\
	\pfe_0^2 \pfe_1 \pfe_0 (u) &=  aq^{-s-2} e_1'^3(u)f_1 + aq^{-s+t} (1 + q^2) e_1'^2(u)\,, \\
	\pfe_0 \pfe_1 \pfe_0^2 (u) &=  aq^{-s-4} e_1'^3(u)f_1 + aq^{-s+t} e_1'^2(u)\,, \\
	\pfe_1 \pfe_0^3 (u) &= aq^{-s-6} e_1'^3(u) f_1\,.
\end{split}
\end{equation*}
Thus, we have $S_{1,0}(u) = 0$ in this case.
\vskip 2mm

By {\it Case 1} and {\it Case 2},  we have seen that $\pfe_0$ and $\pfe_1$ satisfy the quantum Serre relation.
On the other hand,
it is straightforward to check that $\pfe_i$ and $\pfk_j$ satisfy the relation $\pfk_j \pfe_i = q^{a_{ij}} \pfe_i \pfk_j$ for $i,\,j \in \{ 0,\, 1\}$.
Hence the ${\bf k}$-linear operators $\pfe_i$ and $\pfk_i$ $(i = 0, 1)$ give a representation of $U_q({\mf b})$ on $U_q^-({\mf g})$ in the case of $A_1^{(1)}$.
}
\end{ex}

In general, we have the following theorem.

\begin{thm} \label{thm: affine b}
The operators $\mathsf{e}_i$ and $\mathsf{k}_i$ for $0\le i\le n$ satisfy the defining relations of $U_q(\mf{b})$. This gives a representation $\rho_{r,a}^+ : U_q(\mf{b}) \longrightarrow {\rm End}_{\bf k}(U_q^-(\mf{g}))$ given by $\rho_{r,a}^+(e_i)=\mathsf{e}_i$ and $\rho_{r,a}^+(k_i)=\mathsf{k}_i$ for $0\le i\le n$.
\end{thm}
\pf 
We verify that $\mathsf{e}_i$ and $\mathsf{k}_i$ satisfy the defining relations of $U_q(\mf{b})$. 
Note that for type $A_1^{(1)}$, it is done by Example \ref{ex:A1(1)}. Let us consider the remaining cases.
\vskip 1mm

First, we prove the relations for $\pfe_i$ and $\pfk_j$ for $0 \le i, j \le n$.
Clearly, we have $\pfk_i \pfk_j = \pfk_j \pfk_i$ for $0 \le i, j \le n$ by definition.
Let us consider the relations between $\pfe_i$ and  $\pfk_j$ for $0 \le i, j \le n$.
Since $\delta = \alpha_0 + \theta$, we have
$
	(\alpha_0, \alpha_i) = (\delta-\theta, \alpha_i) = -(\theta, \alpha_i).
$
Also, one can check that $(\theta, \theta) = 2$. These implies that 
$
	\pfk_j \pfe_i = q^{a_{ij}} \pfe_i \pfk_j
$
for $0 \le i,\, j \le n$.
\vskip 2mm

Next, we prove the quantum Serre relations of $\mathsf{e}_i$ ($0\le i \le n)$.
By definition, 
$\big\{\, \pfe_i \,\, \big| \,\, i \in I \, \big\}$ satisfies the relations by \cite[Lemma 3.4.2]{Kas91}.
So it remains to show the quantum Serre relations of $\pfe_i$ and $\pfe_0$ for $0\le i \le n$.
Let us consider two cases.
\vskip 1mm

{\em Case 1}. If $i \in J_0$, then by Lemma \ref{lem:dual PBW properties 2}(1), for $u\in U_q^-(\mf{g})_\beta$, we have
\begin{equation*}
	\quad \,\,\,
	\begin{split}
		\pfe_i \pfe_0 (u)
		&= \pfe_i \left( aq^{-( \theta,\, \beta )} u\cdot \xz \right) \\
		&= aq^{-( \theta,\, \beta )} \left\{ e_i' (u) \cdot \xz + q^{( \beta,\,\alpha_i )} u\cdot e_i'(\xz) \right\} \\
		&= aq^{-( \theta,\, \beta + \alpha_i )} \pfe_i (u) \cdot \xz \\
		&= \pfe_0 \pfe_i (u),
	\end{split}	
\end{equation*}
Hence $\pfe_i \pfe_0 = \pfe_0 \pfe_i$ holds.\vskip 1mm

{\em Case 2}. Suppose that $i \in J_1$. 
We claim that
\begin{align} 
	\pfe_0 \pfe_i^2 - (q+q^{-1})\pfe_i \pfe_0 \pfe_i + \pfe_i^2 \pfe_0 = 0,\label{eq:Serre relation 1}\\
	\pfe_i \pfe_0^2 - (q+q^{-1})\pfe_0 \pfe_i \pfe_0 + \pfe_0^2 \pfe_i = 0.\label{eq:Serre relation 2}
\end{align}
Let $u\in U_q^-(\mf{g})_\beta$ be given.
Put 
$s = (\theta,\beta)$ and
$t = (\beta,\alpha_i)$.

First, we have by using Lemma \ref{lem:dual PBW properties 2}
\begin{equation*}
\begin{split}
&\pfe_0 \pfe_i^2 (u) = aq^{-s-2} \left( e_i' \right)^2 (u) \cdot \xz,\\
&\pfe_i \pfe_0 \pfe_i(u)=aq^{-s-1} \left( e_i' \right)^2(u) \cdot \xz + aq^{-s+t+1} e_i'(u) \cdot e_i'(\xz), \\
&\pfe_i^2 \pfe_0 (u)= aq^{-s} \left( e_i' \right)^2(u)\cdot \xz + aq^{-s+t+2} e_i'(u) \cdot e_i' (\xz)+ aq^{-s+t} e_i'(u) \cdot e_i'(\xz).
\end{split}
\end{equation*}
This implies \eqref{eq:Serre relation 1}. Similarly, we have
\begin{equation*}
\begin{split}
&\pfe_i \pfe_0^2 (u) 
			=
			a^2q^{-2s+2} \,e_i'(u)\cdot \xz^2
			+ a^2q^{-2s+t+2}\, u \cdot e_i'(\xz) \cdot \xz
			+ a^2q^{-2s+t} \, u \cdot e_i'(\xz) \cdot \xz,\\
&\pfe_0 \pfe_i \pfe_0 (u)
			=
			a^2q^{-2s+1} e_i'(u)\cdot\xz^2 + a^2q^{-2s+t+1}\, u\cdot e_i'(\xz)\cdot\xz, \\
&\pfe_0^2 \pfe_i (u)
			=
			a^2 q^{-2s} e_i'(u) \cdot \xz^2.
\end{split}
\end{equation*}
This implies \eqref{eq:Serre relation 2}.
\qed\vskip 2mm

\begin{rem} \label{rem:boundary in An(1)}
{\em 
For type $A_n^{(1)}$ with $r = 1$ or $r = n$, since $e_{ \delta_{1,r}\cdot 1 + \delta_{n,r}\cdot n}' (\xz)  = 0$, for $u \in U_q^-({\mf g})_\beta$, we have
\begin{equation*}
	\pfe_0 \pfe_r (u) = aq^{-(\theta, \beta)-1} e_r'(u) \xz = q^{-1} \pfe_r \pfe_0 (u),
\end{equation*}
which also implies the relations \eqref{eq:Serre relation 1} and \eqref{eq:Serre relation 2}. 
}
\end{rem}
}

Now, we are ready to prove Theorem \ref{thm:main-1}.
\vskip 2mm

\noindent {\em Proof of Theorem \ref{thm:main-1}}.
By definition, we have $\xz \in U_q^-(\varpi_r)$. Hence $U_q^-(\varpi_r)$ is invariant under $\pfe_0$.
By Theorem \ref{thm: affine b}, it is enough to verify that $U_q^-(\varpi_r)$ is invariant under $e_i'$ for $i\in I$.
Suppose that $e_i'F^{\rm up}(\beta_k)\neq 0$. 
Note that $ \vep_i ({\bf 1}_{\beta_k})=1 $ by \eqref{eq:epsilon i on beta}.
By Lemma \ref{lem:dual PBW properties 1},
we have
\begin{equation}\label{eq:invariance of root vector}
e_i'F^{\rm up}(\beta_k)=F^{\rm up}(\beta_k - \alpha_i)\in U_q^-(\varpi_r).
\end{equation}
Since $e_i'$ is the derivation on $U_q^-(\mf{g})$ and $U_q^-(\varpi_r)$ is the subalgebra generated by $F^{\rm up}(\beta_k)$, $U_q^-(\varpi_r)$ is invariant under $e_i'$ by \eqref{eq:invariance of root vector}. 
Hence $U_q^-(\varpi_r)$ is a $U_q(\mf{b})$-module.

On the other hand, by definition \eqref{eq:definition of Uwr}, $U_q^-(\varpi_r)$ is ${\mf t}$-diagonalizable and each weight space is finite-dimensional. 
Also, all weights of $U_q^-(\varpi_r)$ are in $D(\overline{\bf 1})$, where $\overline{\bf 1}$ is the map from $I$ to ${\bf k}^\times$ by $i \mapsto 1$ for all $i \in I$.
Hence, $U_q^-(\varpi_r)$ satisfies the conditions in Definition  \ref{df:category O}, that is, $U_q^-(\varpi_r)$ is a $U_q(\mf{b})$-module in the category $\mc{O}$.
\qed

\section{Realization of prefundamental representations} \label{sec:realization}

\subsection{Affine PBW vectors and Drinfeld generators, $\psi_{i, k}^+$}

Let $2\rho$ be the sum of positive roots of $\mathring{\mf g}$.
Take a reduced expression $t_{2\rho}=s_{i_1} \cdots s_{i_N}\in W$ and define a doubly infinite sequence
\begin{equation} \label{eq:doubly infinite sequence}
	\dots, i_{-2}, \,i_{-1}, \,i_0, \,i_1, \, i_2 \dots,
\end{equation}
by setting $i_k = i_{k\, ({\rm mod}\, N)}$ for $k \in \mathbb{Z}$.
Then \eqref{eq:doubly infinite sequence} gives a convex order $\prec$ (cf.~\cite[Definition 2.1]{MT}) on the set of positive roots of $\mf{g}$;
\begin{equation}\label{eq:convex order}
\beta_0 \prec \beta_{-1} \prec \dots \prec \delta \prec \cdots \prec \beta_2 \prec \beta_1, 
\end{equation}
where $\beta_k$ is given by
\begin{equation*}
	\beta_k = 
	\begin{cases} 
	s_{i_0}s_{i_{-1}}\cdots s_{i_{k+1}} (\alpha_{i_k}) & \textrm{if $k \leqq 0$}, \\
	s_{i_1}s_{i_2}\cdots s_{i_{k-1}}(\alpha_{i_k}) & \textrm{if $k > 0$}.
	\end{cases}
\end{equation*}

For $k\in\Z$, we define the root vector $\E_{\beta_k}$ by
\begin{equation} \label{eq:real root vectors}
\E_{\beta_k} = 
\left\{ 
\begin{array}{ll}
T_{i_0}^{-1} T_{i_{-1}}^{-1} \cdots T_{i_{k+1}}^{-1}(e_{i_k}) & \textrm{if $k \leqq 0$}, \\
T_{i_1} T_{i_2} \cdots T_{i_{k-1}}(e_{i_k}) & \textrm{if $k > 0$}.
\end{array}
\right.
\end{equation}
We define the root vector $\F_{\beta_k}$ in the same way with $e_{i_k}$ replaced by $f_{i_k}$ in \eqref{eq:real root vectors}.
Recall that $\E_{\beta_k} \in U_q^+ (\g)$ and $\F_{\beta_k} \in U_q^- (\g)$ for $k \in \mathbb{Z}$ by \cite[Proposition 40.1.3]{Lu93}.
In particular, if $\beta_k = \alpha_i$ for $0\le i \le n$, then $\E_{\beta_k} = e_i$ and $\F_{\beta_k} = f_i$ (cf. \cite[Corollary 4.3]{MT}).

\begin{ex} \label{ex:affine PBW}
	{\rm 
	Let us consider the case of type $A_3^{(1)}$.
	By Remark \ref{rem:formula for w_i}, we have
	\begin{equation*}
	\begin{split}
		t_{\varpi_1} &= \tau_1 (s_3s_2s_1), \\
		t_{\varpi_2} &= \tau_2 (s_2s_1s_3s_2), \\
		t_{\varpi_3} &= \tau_3 (s_1s_2s_3),
	\end{split}
	\end{equation*}
	where $\tau_i$ is the diagram automorphism by the correspondence $\alpha_k \, \mapsto \, \alpha_{k + i\, \text{(mod $4$)}}$.
	Recall that $\rho = \sum_{i=1}^3 \varpi_i$.
	Then the reduced expression $(i_1, \dots, i_{20})$ of $t_{2\rho}$ associated to the following expression (cf. \cite{BCP})
	\begin{equation} \label{eq:exp of t 2rho}
	t_{2\rho} = t_{\varpi_1}t_{\varpi_2}t_{\varpi_3}t_{\varpi_1}t_{\varpi_2}t_{\varpi_3}
	\end{equation}
	is given by 
	\begin{equation*}
	\begin{split}
		(i_1,\, i_2,\, i_3,\, i_4,\, i_5,\, i_6,\, i_7,\, i_8,\, i_9,\, i_{10}) &= (0, 3, {\bf 2}, 1, 0, 2, {\bf 1}, 3, 0, {\bf 1}), \\
		(i_{11},\, i_{12}.\, i_{13},\, i_{14},\, i_{15},\, i_{16},\, i_{17},\, i_{18},\, i_{19},\, i_{20}) &= (2, 1, 0, 3, 2, 0, 3, 1, 2, 3).
	\end{split}
	\end{equation*}
	
	Consider the doubly infinite sequence \eqref{eq:doubly infinite sequence} induced from the above reduced expression.
	With respect to the corresponding convex order $\prec$ \eqref{eq:convex order},  the positive real roots of the form $\delta-\alpha_i$ $(i=1,2,3)$ occur at the positions of the above bold numbers with 
	\begin{equation*}
	\cdots \, \prec \, \delta-\alpha_3 \, \prec \, \cdots \, \prec \, \delta-\alpha_2 \, \prec \, \cdots \, \prec \, \delta-\alpha_1 \, \prec \, \cdots.
	\end{equation*}
	In particular, the root vectors $\E_{\delta-\alpha_i}$ $(i = 1, 2, 3)$ are given by
	\begin{equation} \label{eq:formula of E in A3}
	\begin{split}
		\E_{\delta-\alpha_1} &= e_0 e_3 e_2 - q^{-1} e_3 e_0 e_2 - q^{-1} e_2 e_0 e_3 + q^{-2} e_2 e_3 e_0, \\
		\E_{\delta-\alpha_2} &= e_0 e_1 e_3 - q^{-1} e_1 e_0 e_3 - q^{-1} e_3 e_0 e_1 + q^{-2} e_3 e_1 e_0, \\
		\E_{\delta-\alpha_3} &= e_0 e_1 e_2 - q^{-1} e_1 e_0 e_2 - q^{-1} e_2 e_0 e_1 + q^{-2} e_2 e_1 e_0.
	\end{split}
	\end{equation}
	
	One can observe that the above formula \eqref{eq:formula of E in A3} of $\E_{\delta-\alpha_i}$ does not depend on the order of multiplication in an expression of $t_{2\rho}$ by $t_{\varpi_i}$ (such as \eqref{eq:exp of t 2rho}). 
	For example, let us consider the following expression
\begin{equation*}
	t_{2\rho} = t_{\varpi_2}t_{\varpi_1}t_{\varpi_3}t_{\varpi_2}t_{\varpi_1}t_{\varpi_3},
\end{equation*}
and then the reduced expression $(j_1, \dots, j_{20})$ of $t_{2\rho}$ is given by
\begin{equation*}
\begin{split}
		(j_1,\, j_2,\, j_3,\, j_4,\, j_5,\, j_6,\, j_7,\, j_8,\, j_9,\, j_{10}) &= (0, 3, 1, {\bf 0}, 2, 1, {\bf 0}, 3, 0, {\bf 1}), \\
		(j_{11},\, j_{12}.\, j_{13},\, j_{14},\, j_{15},\, j_{16},\, j_{17},\, j_{18},\, j_{19},\, j_{20}) &= (2, 1, 3, 2, 0, 3, 2, 1, 2, 3).
	\end{split}
\end{equation*}
Here the positive real roots of the form $\delta-\alpha_i$ $(i = 1, 2, 3)$ occur at the positions of the above bold numbers. With respect to the corresponding convex ordering $\prec'$, we have
\begin{equation*}
		\cdots \, \prec' \, \delta-\alpha_3 \, \prec' \, \cdots \, \prec' \, \delta-\alpha_1 \, \prec' \, \cdots \, \prec' \, \delta-\alpha_2 \, \prec' \, \cdots.
\end{equation*}
If we take the doubly infinite sequence \eqref{eq:doubly infinite sequence} associated to $(j_1, \dots, j_{20})$, then by a direct computation, one check that the formulas of $\E_{\delta-\alpha_i}$ are equal to the ones in \eqref{eq:formula of E in A3}.
	}
\end{ex}

In general, we have the following.

\begin{lem} \label{lem:independence of E}
	For $i \in I$ and $k > 0$, the root vectors $\E_{k\delta\pm\alpha_i}$ and $\F_{k\delta\pm\alpha_i}$ are independent of the choice of a reduced expression of $t_{2\rho}$.
\end{lem}
\pf
Let $\prec$ be a convex order associated with $t_{2\rho}$ as above.
One can check that $\delta \prec \beta$ if and only if $\beta=k\delta-\alpha$ for some $k \ge 1$ and a positive root $\alpha$ of $\mathring{\g}$ (cf.~\cite[Lemma 1.1]{BCP}), that is, the convex order $\prec$ has coarse type $w_0$ (see \cite[Definition 2.9]{MT}).
Then the assertion follows from \cite[Corollary 4.6]{MT}.
\qed
\vskip 2mm

Let $o : I \,\rightarrow\, \{ \pm 1 \}$ be a map such that $o(i) = -o(j)$ whenever $a_{ij} < 0$.
Recall that
\begin{equation}\label{eq:Drinfeld generators in terms of affine root vectors}
	\psi_{i, k}^+ = o(i)^k (q-q^{-1})C^{+ \frac{k}{2}}k_i
	\left( \E_{k\delta-\alpha_i}\E_{\alpha_i} - q^{-2}\E_{\alpha_i} \E_{k\delta-\alpha_i} \right),
\end{equation}
for $i \in I$ and $k > 0$ (see \cite[Proposition 1.2]{BCP}).
We remark that the right-hand side of \eqref{eq:Drinfeld generators in terms of affine root vectors} is also independent of the choice of a reduced expression of $t_{2\rho}$ due to Lemma \ref{lem:independence of E}.

The following lemma together with \eqref{eq:Drinfeld generators in terms of affine root vectors} enables us to compute the action of $\psi_{i, k}^+$ for $i \in I$ and $k \in \mathbb{Z}_{>0}$ on a $U_q({\mf b})$-module (cf. Example \ref{ex:l-weight in A1}, \ref{ex:l-weight in A3} and \ref{ex:l-weight in D4}).

\begin{lem} \label{lem:inductive formula}
For $i \in I$ and $k \in \mathbb{Z}_{>0}$, we have
\begin{equation*}
\begin{split}
\E_{(k+1)\delta-\alpha_i} &= -\frac{1}{q+q^{-1}}
\Big(\,
\E_{\delta-\alpha_i} \E_{\alpha_i}\E_{k\delta-\alpha_i} - q^{-2} \E_{\alpha_i}\E_{\delta-\alpha_i}\E_{k\delta-\alpha_i} \\
& \qquad \qquad \quad\,\,\,- \E_{k\delta -\alpha_i} \E_{\delta-\alpha_i} \E_{\alpha_i} + q^{-2} \E_{k\delta - \alpha_i} \E_{\alpha_i} \E_{\delta - \alpha_i}\,
\Big)\,.
\end{split}
\end{equation*}
\end{lem}
\pf Thanks to Lemma \ref{lem:independence of E}, we may use the following relation \cite[Proposition 1.2]{BCP} (cf. \cite{B}): for $i,\,j \in I$ and $k,\, l \in \mathbb{Z}_+$ with $l > 0$,
\begin{equation} \label{eq:inductive formula for root vectors}
\left[ \E_{l\delta, i}, \, \E_{k\delta\pm\alpha_j} \right] = \pm o(i)^l o(j)^l \frac{1}{l} \left[ la_{ij} \right]_q \E_{(l+k)\delta \pm \alpha_i},
\end{equation}
where $k > 0$ for $k\delta-\alpha_j$. Here $\E_{l\delta, i} \in U_q(\g)$ is defined by the functional equation
\begin{equation} \label{eq:functional equation}
\exp\left( (q-q^{-1}) \sum_{l=1}^{\infty} \E_{l\delta, i} u^l \right) = 1 + \sum_{k=1}^{\infty} (q-q^{-1}) \wtd{\psi}_{i, k} u^k\,,
\end{equation}
where $\wtd{\psi}_{i, k} = \E_{k\delta-\alpha_i}\E_{\alpha_i} - q^{-2}\E_{\alpha_i} \E_{k\delta-\alpha_i}$.
In particular, by comparing the coefficients of $u$ in \eqref{eq:functional equation}, we get
\begin{equation} \label{eq:E delta i}
\E_{\delta, i} = \wtd{\psi}_{i, 1} = \E_{\delta-\alpha_i} \E_{\alpha_i} - q^{-2} \E_{\alpha_i} \E_{\delta-\alpha_i},
\end{equation}
On the other hand, if we put $l = 1$ and $i = j$ in \eqref{eq:inductive formula for root vectors}, then for $k > 0$, we have
\begin{equation} \label{eq:inductive formula in low rank}
-[2] \E_{(k+1)\delta-\alpha_i} 
= 
\E_{\delta, i}\E_{k\delta-\alpha_i} - \E_{k\delta-\alpha_i}\E_{\delta, i},
\end{equation}
which implies the required formula by \eqref{eq:E delta i} and \eqref{eq:inductive formula in low rank}.
\qed\vskip 2mm

\subsection{Highest $\ell$-weight of $U_q^-(\varpi_r)_a$}\label{subsec:loop highest weight}
Let $U_q^-(\varpi_r)_a$ denote the $U_q(\mf{b})$-module defined by \eqref{eq:borel action} in Theorem \ref{thm: affine b}. We denote by $\V \in U_q^-(\varpi_r)_a$ a non-zero vector of weight $0$, which is equal to $1 \in U_q^-(\varpi_r)$ up to a scalar multiplication.
In this subsection, we compute the $\ell$-weight of $\V$ in $U_q^-(\varpi_r)_a$ explicitly.
It is crucial to prove that the $U_q(\mf{b})$-module $U_q^-(\varpi_r)_a$ is isomorphic to the prefundamental representation $L_{r,1}^+$ for $r$ in \eqref{eq:minuscule index} up to a shift of spectral parameter.

\subsubsection{\bf Type $A_n^{(1)}$} \label{subsec:l-weight in A}
Let us first give explicit examples to compute the $\ell$-weight of $\V$. Then the rest of this subsection will be devoted to generalize the following examples to an arbitrary rank $n$ and $r \in I$.

\begin{ex} \label{ex:l-weight in A1}
{\em 
Let us consider the case of type $A_1^{(1)}$. We continue to use the convention in Example \ref{ex:A1(1)}.
By \eqref{eq:Drinfeld generators in terms of affine root vectors}, it suffice to compute the action of $\E_{k\delta-\alpha_1}$ on $\V$ in order to compute the eigenvalue of $\psi_{1, k}^+$ for $\V$.
In this case, $t_{2\rho} = s_0 s_1$. Thus, we have
\begin{equation*}
	\E_{\delta-\alpha_1} = e_0 \quad \text{and} \quad \E_{\delta-\alpha_1} \V = a\cdot \xz = a\cdot f.
\end{equation*}
By Lemma \ref{lem:inductive formula}, the action of $\E_{2\delta-\alpha_1}$ on $\V$ is given by 
\begin{equation*}
	\E_{2\delta-\alpha_1} \V = -\frac{1}{q+q^{-1}} a^2 \left( 1 - (1+q^{-2}) + q^{-2} \right) \cdot f = 0,
\end{equation*}
which implies that $\E_{k\delta-\alpha_1} \V = 0$ for $k \ge 2$ by Lemma \ref{lem:inductive formula}. Hence, we have
\begin{equation*}
	\psi_{1, k}^+ \V =
	\begin{cases}
		-o(1)(q-q^{-1})q^{-2} a \V & \text{if $k = 1$,} \\
		0 & \text{if $k > 1$.}
	\end{cases}
\end{equation*}
Finally, the $\ell$-weight $\Psi(z)$ of $\V$ is given by $(\,1-o(1)(q-q^{-1})q^{-2}az\,)$.
}
\end{ex}

\begin{ex} \label{ex:l-weight in A3}
{\em 
As another example, we consider the case of $A_3^{(1)}$ with $r = 2$.
Let us recall Example \ref{ex:affine PBW}.
\vskip 2mm

{\em Step 1}. 
By Lemma \ref{lem:dual PBW properties 1} and \eqref{eq:formula of E in A3},
\begin{equation*}
\begin{split}
	\E_{\delta-\alpha_2}\V &= q^{-2} \pfe_3 \pfe_1 \pfe_0 \V = q^{-2} a \pfe_3 \pfe_1 (\xz) = q^{-2} a (f_2).
\end{split}
\end{equation*}
Also, we have
\begin{equation*}
\begin{split}
	\E_{\delta-\alpha_2} (f_2) &= q^{-2} \pfe_3 \pfe_1 \pfe_0 (f_2) 
	= q^{-2} a \pfe_3 \pfe_1 (f_2 \xz) 
	= a (f_2 e_3'e_1'(\xz)) = a (f_2^2).
\end{split}
\end{equation*}
\vskip 2mm

{\em Step 2.}
By Lemma \ref{lem:inductive formula} and the above computation, we observe that
\begin{equation*}
	\E_{2\delta-\alpha_2} \V = -\frac{1}{q+q^{-1}} a^2 \left( q^{-4} - q^{-4} (1 + q^{-2} ) + q^{-6}  \right) f_2 = 0.
\end{equation*}
Then it follows from Lemma \ref{lem:inductive formula} that
$\E_{k\delta-\alpha_2} \V = 0$
for all $k \ge 2$.
\vskip 2mm

{\em Step 3.}
By \eqref{eq:Drinfeld generators in terms of affine root vectors}, the eigenvalue of $\V$ under $\psi_{2, k}^+$ is given by 
\begin{equation*}
	\psi_{2, k}^+ \V =
	\begin{cases}
		-o(2)(q-q^{-1})q^{-4}a \V & \text{if $k = 1$,} \\
		0 & \text{if $k > 1$.}
	\end{cases}
\end{equation*}

For $i = 1, 3$, we have $\E_{\delta-\alpha_i} \V = 0$ by \eqref{eq:formula of E in A3} and Lemma \ref{lem:dual PBW properties 1}(2). By Lemma \ref{lem:inductive formula}, we conclude $\Psi_1(z) = \Psi_3(z) = 1$. Finally, the $\ell$-weight $\Psi(z)$ of $\V$ is 
$$\Psi(z) = (1,\, 1-o(2)(q-q^{-1})q^{-4}az,\, 1).$$
}
\end{ex}

Now, let us consider the case of type $A_n^{(1)}$ $(n \ge 2)$ with $r \in I$. 
For simplicity, we denote by $\pm q^{-\mathbb{Z}_+}$ the subset of ${\bf k}$ consisting of elements of the form $\pm q^{-m}$ for $m \in \mathbb{Z}_+$.

\begin{lem} \label{Lem: x wr}
Let $s_{i_1} s_{i_2} \cdots s_{i_\ell}$ be a reduced expression of $w_r^{-1}$, where $w_r$ is given in Lemma \ref{lem:property of w_r}.
Then $x := T_{i_1} T_{i_2} \cdots T_{i_{\ell - 1}} (e_r) $ can be written as 
\begin{equation*}
 (-q^{-1})^{n-1} (e_{1}\dots e_{n-r-1}e_{n-r})(e_{n}\dots e_{n+3-r}e_{n+2-r})e_{n+1-r} + \sum_{j_1, \dots, j_n}
 C_{j_1, \dots, j_n}(q)
 e_{j_1}\dots e_{j_n},
\end{equation*}
where the sum is over the sequences $(j_1,\dots,j_n)$ such that $\sum_{k=1}^n \alpha_{j_k} = \delta-\alpha_0$ with $j_n \neq n+1-r$ and 
$C_{j_1, \dots, j_n}(q) \in \pm q^{-\mathbb{Z}_+}$.
\end{lem}
\pf
By Lemma \ref{lem:property of w_r}, we have 
\begin{itemize}
\item $w_r^{-1} = w_{n+1-r}$,
\item  $x$ does not depend on the choice of reduced expressions of $w_r^{-1}$. 
\end{itemize} 
Thus we use the reduced expression in Remark \ref{rem:formula for w_i}. 
Set $[a,b]_q := ab - q^{-1} ba$. Note that $T_i(e_j) = [e_i, e_j]_{q}$ if $a_{ij}=-1$.
We shall use induction on rank $n$. Since it is obvious when $n=1$, we assume that $n>1$. 

Suppose that $r=1$. Since $w_1^{-1} = w_{n} = s_n s_{n-1} \cdots s_1$, we have 
\begin{align*}
x& = T_n \cdots T_4 T_3 T_2 (e_1)	= T_n \cdots T_4 T_3 ( [e_2, e_1]_q)
= T_n \cdots T_4  ( [ [e_3, e_2]_q, e_1]_q) \\
&= [  \cdots  [[ e_n, e_{n-1}]_q, e_{n-2}  ]_q, \cdots e_1 ]_q,
\end{align*}
which implies the assertion.

Suppose that $r > 1$. For a reduced expression $w = s_{j_1} \cdots s_{j_k}$, we let
$$
 w^- = s_{j_1} \cdots s_{j_{k-1}} \qquad \text{and} \qquad  T_{w} = T_{j_1} \cdots T_{j_k}.
 $$ 
 
Let $v = {\bf s}_{(n-r,\, n-1)} {\bf s}_{(n-r-1,\, n-2)} \dots {\bf s}_{(1,\, r)}$ and 
$
y = T_{v^-  } (e_r). 
$ 
Note that $v \in \langle \, s_1, \dots, s_{n-1} \, \rangle $ and 
$ x= T_{ {\bf s}_{(n+1-r,\, n)} } y$.
By the induction hypothesis, $y$ can be written as 
	\begin{equation*}
	(-q^{-1})^{n-2} (e_{1}\dots e_{n-r-2}e_{n-r-1})(e_{n-1}\dots e_{n+2-r}e_{n+1-r})e_{n-r} + A,
\end{equation*}
where $A$ is of the form $\sum_{j_1, \dots, j_{n-1}}C_{j_1, \dots, j_{n-1}}(q) e_{j_1}\dots e_{j_{n-1}}$  with $j_{n-1} \neq n-r$.
As
\begin{equation} \label{eq:T from n+1-r to n}
T_{ {\bf s}_{(n+1-r,\, n)} } = T_{ n+1-r} T_{n+2-r} \cdots T_n,
\end{equation}
it follows from \eqref{eq:T from n+1-r to n} that 
\begin{align} \label{eq:action of Tw on e}
	T_{ {\bf s}_{(n+1-r,\, n)} } ( e_k) = 
	\begin{cases}
		e_k & \text{ if } k=1, \ldots, n-r-1, \\
		[e_{n+1-r}, e_{n-r}]_q & \text{ if } k=  n-r, \\
		e_{k+1} & \text{ if } k= n+1-r, \ldots, n-1. \\
	\end{cases}
\end{align}
Therefore, we have 
\begin{align*}
x &= T_{ {\bf s}_{(n+1-r,\, n)} } ( y)	\\
&= 	(-q^{-1})^{n-2} (e_{1}\dots e_{n-r-2}e_{n-r-1})(e_{n}\dots e_{n+3-r}e_{n+2-r})	[e_{n+1-r}, e_{n-r}]_q  + T_{ {\bf s}_{(n+1-r,\, n)} } (A) \\
&= 	(-q^{-1})^{n-1} (e_{1}\dots e_{n-r-2}e_{n-r-1}  )(e_{n}\dots e_{n+3-r}e_{n+2-r})	e_{n-r} e_{n+1-r}  \\
& \quad  + (-q^{-1})^{n-2} (e_{1}\dots e_{n-r-2}e_{n-r-1})(e_{n}\dots e_{n+3-r}e_{n+2-r})	e_{n+1-r} e_{n-r} + T_{ {\bf s}_{(n+1-r,\, n)} } (A),
\end{align*}
which implies the assertion since the last two terms satisfy the required property by \eqref{eq:action of Tw on e}.\!\!\!
\qed

The following formula plays an important role in the computation of the $\ell$-weight of $\V\in U_q^-(\varpi_r)_a$.
\vskip 1mm

\begin{lem} \label{lem:formula of E}
For $r \in I$, we have
	\begin{equation*}
		\E_{\delta-\alpha_r} = (-q^{-1})^{n-1} (e_{r+1}\dots e_{n-1}e_n)(e_{r-1}\dots e_2e_1)e_0 + \sum_{j_1, \dots, j_n} C_{j_1, \dots, j_n}(q) e_{j_1}\dots e_{j_n},
	\end{equation*}
where the sum is over the sequences $(j_1,\dots,j_n)$ such that $\sum_{k=1}^n \alpha_{j_k} = \delta-\alpha_r$ with $j_n \neq 0$ and $C_{j_1, \dots, j_n}(q) \in \pm q^{-\mathbb{Z}_+}$.
\end{lem}
\pf
Thanks to Lemma \ref{lem:independence of E}, $\E_{\delta-\alpha_r}$ does not depend on the choice of a reduced expression of $t_{2\rho}$. By Lemma \ref{lem:property of w_r}, we can write 
$$
t_{2\rho} = t_{ \varpi_{ r }} t_{2\rho - \varpi_r} = \tau_{r} w_{r}^{-1} t_{2\rho - \varpi_r},
$$
where $\tau_r$ is the Dynkin diagram automorphism of $A_n^{(1)}$ sending $i$ to $i+r \quad  ( \mathrm{mod}\  n+1)$.
Let $s_{i_1} s_{i_2} \cdots s_{i_\ell}$ be a reduced expression of $w_r^{-1}$. Note that $w_r^{-1} = w_{n+1-r}$ and $ i_\ell = r$ by Lemma \ref{lem:property of w_r}.  Since $ t_{\varpi_r}(\alpha_r) =  \alpha_r - \delta$, we know that 
$$
\E_{\delta-\alpha_r} = \tau_{r}  T_{i_1} T_{i_2} \cdots T_{i_{\ell - 1}} (e_r).
$$
Here we understand $\tau_r$ as the automorphism of $U_q(\g)$ induced by the corresponding Dynkin diagram automorphism. 
Thus the assertion follows from Lemma \ref{Lem: x wr}.
\qed

\begin{cor} \label{cor:action of real root vector}	
	The action of $\E_{\delta-\alpha_i}$ on ${\bf 1} \in U_q^-(\varpi_r)_a$ is given by
	\begin{equation*}
		\E_{\delta-\alpha_i} {\bf 1} = 
		\begin{cases}
			(-q^{-1})^{n-1} a \left( f_r \right)  & \text{if $i = r$}, \\
			0 & \text{if $i \neq r$}.
		\end{cases}
	\end{equation*}
	Moreover, we have $\E_{k\delta-\alpha_i} {\bf 1} = 0$ for $i \in I$ and $k \ge 2$.
\end{cor}
\pf
We use Lemma \ref{lem:dual PBW properties 1} and Lemma \ref{lem:formula of E} to compute the action of $\E_{\delta-\alpha_i} {\bf 1}$. 
For $i = r$, we have by Lemma \ref{lem:formula of E},
\begin{equation*}
	\E_{\delta-\alpha_r} \V = (-q^{-1})^{n-1} (\mathsf{e}_{r+1} \dots \mathsf{e}_{n-1} \mathsf{e}_n) (\mathsf{e}_{r-1} \dots \mathsf{e}_2 \mathsf{e}_1) \mathsf{e}_0 \V.
\end{equation*}
By Lemma \ref{lem:dual PBW properties 1} and \eqref{eq:borel action}, 
it is enough to consider the transition of the weight of $\xz$ along the action $(\mathsf{e}_{r+1} \dots \mathsf{e}_{n-1} \mathsf{e}_n) (\mathsf{e}_{r-1} \dots \mathsf{e}_2 \mathsf{e}_1)$. Consequently, we have $\E_{\delta-\alpha_r} \V = (-q^{-1})^{n-1} a \left( f_r  \right)$.
For $i \neq r$, since $-\alpha_i$ is not a weight of $U_q^-(\varpi_r)$, we conclude that $\E_{\delta-\alpha_i} \V = 0$ for $i \neq r$.
\qed

\begin{rem}
{\em 
In the proof of Corollary \ref{cor:action of real root vector}, for $i \neq r$, one can check $\E_{\delta-\alpha_i} \V = 0$ directly by using Lemma \ref{lem:formula of E} since it holds for any $r \in I$.
}
\end{rem}

\vskip 2mm

\begin{cor} \label{cor:action of f_r}
The action of $\E_{\delta-\alpha_r}$ on $f_r \in U_q^-(\varpi_r)_a$ is given by
	\begin{equation*}
		\E_{\delta-\alpha_r} (f_r) = (-q^{-1})^{n-1} q^2 a \left( f_r^2 \right).
	\end{equation*}
\end{cor}
\pf
We use Lemma \ref{lem:formula of E} to compute the action of $\E_{\delta-\alpha_i}$ on $f_r$.
Since $e_i' f_r = 0$ for $i \neq r$, the summation of $e_{j_1} \dots e_{j_n}$ $(j_n \neq 0)$ in the formula of $\E_{\delta-\alpha_r}$ in Lemma \ref{lem:formula of E} is vanished when acting on $f_r$. Thus, we have
\begin{equation*}
	\E_{\delta-\alpha_r} (f_r) = (-q^{-1})^{n-1} (\mathsf{e}_{r+1} \dots \mathsf{e}_{n-1} \mathsf{e}_n) (\mathsf{e}_{r-1} \dots \mathsf{e}_2 \mathsf{e}_1) \mathsf{e}_0 (f_r).
\end{equation*}
Then since the action of $\pfe_i$ $(i \neq r)$ is defined by the derivation \eqref{eq:derivation}, we have 
\begin{equation*}
	\pfe_i (f_r F^{\rm up}(\beta)) = q^{-(\alpha_r, \alpha_i)} f_r \left( e_i' F^{\rm up}(\beta) \right).
\end{equation*}
The assertion follows from Lemma \ref{lem:dual PBW properties 1}.
\qed

\begin{prop}\label{prop:highest l-weight of type A}
	The $\ell$-weight $\Psi=(\Psi_i(z))_{i\in I}$ of $\V \in U_q^-(\varpi_r)_a$ is given by
	\begin{equation*} %\label{eq:l-weight}
		\Psi_i(z) = 
		\begin{cases}
			1-ac_{r}z & i = r, \\
			1 & i \neq r,
		\end{cases}
	\end{equation*}
	where $c_{r} = (-1)^{n+1} o(r) (q-q^{-1}) q^{-n-1}$.
\end{prop}
\pf
By \eqref{eq:Drinfeld generators in terms of affine root vectors} and Corollary \ref{cor:action of real root vector},
\begin{equation*}
	\psi_{i, 1}^+ {\bf 1} =
	\begin{cases}
		(-1)^n o(r)(q-q^{-1})q^{-n-1} a {\bf 1} & \text{if $i = r$}, \\
		0 & \text{if $i \neq r$}.
	\end{cases}
\end{equation*}

We claim that $\psi_{i, k}^+ {\bf 1} = 0$ for all $k \ge 2$. It suffice to show that $\E_{k\delta-\alpha_i} \V = 0$ for all $k \ge 2$ due to the relation \eqref{eq:Drinfeld generators in terms of affine root vectors}. Then the proof is by induction on $k$.
When $k=2$, by Lemma \ref{lem:inductive formula}, we have
\begin{equation*}
\begin{split}
\E_{2\delta-\alpha_i} &= -\frac{1}{q+q^{-1}}
\Big(\,
\E_{\delta-\alpha_i} \E_{\alpha_i}\E_{\delta-\alpha_i} - q^{-2} \E_{\alpha_i}\E_{\delta-\alpha_i}\E_{\delta-\alpha_i} \\
& \qquad \qquad \quad\,\,\,- \E_{\delta -\alpha_i} \E_{\delta-\alpha_i} \E_{\alpha_i} + q^{-2} \E_{\delta - \alpha_i} \E_{\alpha_i} \E_{\delta - \alpha_i}\,
\Big)\,.
\end{split}
\end{equation*}
If $i \neq r$, we have $\E_{2\delta-\alpha_i} \V = 0$ by Corollary \ref{cor:action of real root vector}.
Suppose that $i = r$. By the above formula, we have
\begin{equation} \label{eq:inductive formula for k=2}
	\E_{2\delta-\alpha_r} \V = 
	-\frac{1}{q+q^{-1}}
\Big(\,
\E_{\delta-\alpha_r} \E_{\alpha_r}\E_{\delta-\alpha_r} - q^{-2} \E_{\alpha_r}\E_{\delta-\alpha_r}\E_{\delta-\alpha_r} + q^{-2} \E_{\delta - \alpha_r} \E_{\alpha_r} \E_{\delta - \alpha_r}\,
\Big)\V.
\end{equation}
By Corollaries \ref{cor:action of real root vector} and \ref{cor:action of f_r},
\begin{equation*}
\begin{split}
	\E_{\delta-\alpha_r} \E_{\alpha_r}\E_{\delta-\alpha_r}\V &= q^{-2n+2} a^2 (f_r\V), \quad
	\E_{\alpha_r}\E_{\delta-\alpha_r}\E_{\delta-\alpha_r}\V = q^{-2n+4} (1+q^{-2})a^2 (f_r).
\end{split}
\end{equation*}
Putting them in \eqref{eq:inductive formula for k=2}, we conclude that $\E_{2\delta-\alpha_r} \V = 0$.
For $k > 2$, the induction step follows directly from Lemma \ref{lem:inductive formula}. This completes the proof.
\qed

\subsubsection{\bf Type $D_n^{(1)}$} \label{subsec:l-weight in D}
Let us compute the $\ell$-weight of $\V \in U_q^-(\varpi_r)_a$ in the case of type $D_n^{(1)}$ with $r  \in \{\, 1,\, n-1,\, n\,\}$ following the arguments in Section \ref{subsec:l-weight in A}.

\begin{lem} \label{lem:formula of E for type D}
\mbox{} \
\begin{itemize}
	\item[(1)] If $r = 1$, then
		\begin{equation*}
		\begin{split}
			\E_{\delta-\alpha_1} &= q^{-2n+4} (e_2 e_3 \dots e_{n-2} e_{n-1})(e_n e_{n-2} \dots e_3 e_2)e_0 \\
			&\hspace{3.5cm} + \sum_{j_1, j_2, \dots, j_{2n-3}} C_{j_1, j_2, \dots, j_{2n-3}}(q) e_{j_1}e_{j_2}\dots e_{j_{2n-3}},
		\end{split}
		\end{equation*}
		where the sum is over the sequences $(j_1, j_2, \dots,j_{2n-3})$ such that $\sum_{k=1}^{2n-3} \alpha_{j_k} = \delta-\alpha_1$ with $j_{2n-3} \neq 0$ and $C_{j_1, j_2, \dots, j_{2n-3}}(q) \in \pm q^{-\mathbb{Z}_+}$.
	\vskip 2mm
	
	\item[(2)] If $r = n$, then
		\begin{equation*}
\begin{split}
\E_{\delta-\alpha_n} &=\, q^{-2n+4} (e_{n-2}e_{n-3}\dots e_2e_1)(e_{n-1}e_{n-2}\dots e_3e_2)e_0 \\ 
&\hspace{3.5cm} + \sum_{j_1, j_2, \dots, j_{2n-3}} C_{j_1, j_2, \dots, j_{2n-3}}(q) e_{j_1}e_{j_2}\dots e_{j_{2n-3}},
\end{split}
\end{equation*}
where the sum is over the sequences $(j_1, j_2, \dots,j_{2n-3})$ such that $\sum_{k=1}^{2n-3} \alpha_{j_k} = \delta-\alpha_n$ with $j_{2n-3} \neq 0$ and $C_{j_1, j_2, \dots, j_{2n-3}}(q) \in \pm q^{-\mathbb{Z}_+}$.
	\vskip 2mm
	
	\item[(3)] If $r = n-1$, then the formula of $\E_{\delta-\alpha_{n-1}}$ is obtained from the one of $\E_{\delta-\alpha_n}$ by replacing $e_{n-1}$ with $e_n$ in which the sum is over the sequences $(j_1, j_2, \dots,j_{2n-3})$ such that $\sum_{k=1}^{2n-3} \alpha_{j_k} = \delta-\alpha_{n-1}$ with $j_{2n-3} \neq 0$ and $C_{j_1, j_2, \dots, j_{2n-3}}(q) \in \pm q^{-\mathbb{Z}_+}$.
\end{itemize}
\end{lem}
\pf
Let $\mathring{\mf g}'$ be the simple Lie algebra of type $D_{n-1}$ associated to $I \, \setminus \, \{\, 1 \,\}$ and let $\varpi_i'$ be the $i$-th fundamental weight of $\mathring{\mf g}'$.
Put $t_{\varpi_r'}^{-1} = w_r \tau$ for some $\tau \in \mathcal{T}'$.
Note that the reduced expression ${\bf i}$ of $w_r^{-1}$ can be rewritten by
\begin{equation} \label{eq:induction on i in D}
	{\bf i} =
	\begin{cases}
		1 \cdot {\bf i}' \cdot 1 & \text{if $r = 1$,} \\
		(r,\, r-2,\, \cdots,\, 1) \cdot {\bf i}' & \text{if $r=n$ and $r$ is even,} \\
		(r-1,\, r-2\, \cdots,\, 1) \cdot {\bf i}' & \text{if $r=n$ and $r$ is odd,}
	\end{cases}
\end{equation}
where ${\bf i}'$ is the sequence obtained from the reduced expression of $(w_r')^{-1}$ by shifting the index by $1$ (cf. Example \ref{ex:induction on formula of E in D}).
Then the formulas in (1) and (2) follow from the same argument as in Lemma \ref{Lem: x wr} and Lemma \ref{lem:formula of E} by using
the reduced expression of $w_r$ in Remark \ref{rem:formula for w_i} (see also $w^J$ with $J = I \setminus \{ n \}$ in \cite[Section 3.1]{JK19}) and \eqref{eq:induction on i in D}.
Also, the formula in (3) is proved following the previous computation for the formula in (2).
\qed
\vskip 1mm

\begin{ex} \label{ex:root vector in D4}
{\em 
Let us consider the case of $D_4^{(1)}$ with $r = 1$ or $4$.
\vskip 1mm

{\em Case 1.} $r = 1$.
In this case, we have
\begin{equation} \label{eq:w1 in D4}
	t_{\varpi_1} = \tau_1 (s_1 s_2) (s_4 s_3 s_2 s_1),
\end{equation}
where $\tau_1$ is the Dynkin diagram automorphism so that $\alpha_0 \, \leftrightarrow \, \alpha_1$ and $\alpha_3 \, \leftrightarrow \, \alpha_4$ (cf. \cite{Bo}).
By Lemma \ref{lem:independence of E} and \eqref{eq:w1 in D4}, it is easy to check that
\begin{equation*}
	\E_{\delta-\alpha_1} = T_0 T_2 T_3 T_4 T_2 (e_0) = q^{-4}(e_2e_3)(e_4e_2)e_0 + \sum_{j_1, j_2, \dots, j_5} C_{j_1, j_2, \dots, j_5}(q) e_{j_1} e_{j_2} \dots e_{j_5},
\end{equation*}
where the sum is over the sequences $(j_1, j_2, \dots, j_5)$ such that $\sum_{k=1}^{5}a_{j_k} = \delta-\alpha_1$ with $j_5 \neq 0$.
\vskip 1mm

{\em Case 2.} $r = 4$.
In this case, we have
\begin{equation} \label{eq:w4 in D4}
	t_{\varpi_1} = \tau_4 s_4 (s_2 s_3) (s_1 s_2 s_4),
\end{equation}
where $\tau_4$ is the Dynkin diagram automorphism so that $\alpha_0 \, \leftrightarrow \alpha_4$ and $\alpha_1 \, \leftrightarrow \alpha_3$ (cf. \cite{Bo}).
By Lemma \ref{lem:independence of E} and \eqref{eq:w4 in D4}, it is easy to check that
\begin{equation*}
	\E_{\delta-\alpha_4} = T_0 T_2 T_1 T_3 T_2 (e_0) = q^{-4} (e_2 e_1) (e_3 e_2) e_0 + \sum_{j_1, j_2, \dots, j_5} C_{j_1, j_2, \dots, j_5}(q) e_{j_1} e_{j_2} \dots e_{j_5},
\end{equation*}
where the sum is over the sequences $(j_1, j_2, \dots, j_5)$ such that $\sum_{k=1}^{5}a_{j_k} = \delta-\alpha_4$ with $j_5 \neq 0$.
}
\end{ex}

\begin{ex} \label{ex:induction on formula of E in D}
{\em 
Let us consider the case of $D_5^{(1)}$ with $r = 5$. In this case, we have
\begin{equation*}
	t_{\varpi_5} = \tau_5 s_4 (s_3 s_5) (s_2 s_3 s_4) (s_1 s_2 s_3 s_5),
\end{equation*}
where $\tau_5$ is the Dynkin diagram automorphism so that $\alpha_0 \, \leftrightarrow \, \alpha_5$, $\alpha_1 \, \leftrightarrow \, \alpha_4$ and $\alpha_2 \, \leftrightarrow \alpha_3$.
Note that the reduced expression of $\tau_5^{-1} t_{\varpi_5}$ can be rewritten by
\begin{equation*}
	(s_4 s_3 s_2 s_1) (s_5) (s_3 s_4) (s_2 s_3 s_5).
\end{equation*}

By similar computation as in Example \ref{ex:root vector in D4}, one can check that 
\begin{equation*}
	T_5 T_3 T_4 T_2 T_3 (e_5) = q^{-4} (e_3 e_4) (e_2 e_3) e_5 + \sum_{j_1, j_2, \dots, j_5} C_{j_1, j_2, \dots, j_5}(q) e_{j_1} e_{j_2} \dots e_{j_5},
\end{equation*}
where the sum is over the sequences $(j_1, j_2, \dots, j_5)$ such that $\sum_{k=1}^{5}a_{j_k} = \delta-\alpha_1-\alpha_2$ with $j_5 \neq 5$.
Note that $s_4s_3s_2s_1$ sends $\alpha_i$ into $\alpha_{i-1}$ for $i = 2, 3, 4$. By applying $T_4T_3T_2T_1$ to the above formula, we have
\begin{equation*}
	T_{w_5^{-1}s_5}(e_5) = q^{-6} (e_2 e_3 e_5) (e_1 e_2 e_3) (e_4) + \sum_{j_1, j_2, \dots, j_7} C_{j_1, j_2, \dots, j_7}(q) e_{j_1} e_{j_2} \dots e_{j_7},
\end{equation*}
where the sum is over the sequences $(j_1, j_2, \dots, j_7)$ such that $\sum_{k=1}^{7}a_{j_k} = \delta-\alpha_0$ with $j_7 \neq 4$.

Finally, by applying $\tau_5$, we have
\begin{equation*}
\E_{\delta-\alpha_5} = q^{-6} (e_3e_2e_1) (e_4 e_3 e_2) e_0 + \sum_{j_1, \dots, j_7} C_{j_1, j_2, \dots, j_7}(q) e_{j_1} e_{j_2} \dots e_{j_7},
\end{equation*}
where the sum is over the sequences $(j_1, j_2, \dots, j_7)$ such that $\sum_{k=1}^{7}a_{j_k} = \delta-\alpha_5$ with $j_7 \neq 0$.
}
\end{ex}
\vskip 1mm

\begin{cor} \label{cor:action of real root vector for type D}
\
\begin{enumerate}
	\item The action of $\E_{\delta-\alpha_i}$ on $\V \in U_q^-(\varpi_r)_a$ is given as follows:
\begin{equation*}
	\E_{\delta-\alpha_i} \V = 
		\begin{cases}
			aq^{-2n+4} (f_r) & \text{if $i = r$}, \\
			0 & \text{if $i \neq r$}.
		\end{cases}
\end{equation*}
Moreover, we have $\E_{k\delta-\alpha_i} \V = 0$ for $i \in I$ and $k \ge 2$.
\vskip 2mm

	\item 
	The action of $\E_{\delta-\alpha_r}$ on $f_r \in U_q^-(\varpi_r)_a$ is given by
	\begin{equation*}
		\E_{\delta-\alpha_r} (f_r) = aq^{-2n+6} (f_r^2 ).
	\end{equation*}
\end{enumerate}
\end{cor}
\pf
The proofs of (1) and (2) are almost identical with the ones of Corollaries \ref{cor:action of real root vector} and \ref{cor:action of f_r}, respectively,  by using Lemma \ref{lem:formula of E for type D}.
\qed

\begin{prop} \label{prop:highest l-weight of type D}
	The $\ell$-weight $\Psi = (\Psi_i(z))_{i \in I}$ of $\V\in U_q^-(\varpi_r)_a$ is given by
	\begin{equation*} %\label{eq:l-weight}
		\Psi_i(z) = 
		\begin{cases}
			1 - ac_r z & i = r, \\
			1 & i \neq r,
		\end{cases}
	\end{equation*}
	where $c_r = o(r)(q-q^{-1})q^{-2(n-1)}$.
\end{prop}
\pf
By following the proof of Proposition \ref{prop:highest l-weight of type A} with Corollary \ref{cor:action of real root vector for type D},
\begin{equation*}
	\E_{k\delta - \alpha_i} \V = 
	\begin{cases}
		aq^{-2n+4} (f_r) & \text{if $i=r$ and $k=1$,} \\
		0 & \text{otherwise.}
	\end{cases}
\end{equation*}
Then the assertion follows from the relation \eqref{eq:Drinfeld generators in terms of affine root vectors}.
\qed

\subsection{Realization of $L^+_{i,a}$} \label{subsec:realization L+}
Now we are in a position to state the main result in this paper.

\begin{thm}\label{thm:main-2}
For a minuscule $\varpi_r$ and $a\in {\bf k}^\times$, we have 
\begin{equation*}
U_q^-(\varpi_r)_a \cong L_{r,ac_r}^+,
\end{equation*} 
as a $U_q(\mf{b})$-module, where $c_{r}$ is given by
\begin{equation} \label{eq:c_r}
	c_r =
	\begin{cases}
		 (-1)^{n+1} o(r) (q-q^{-1}) q^{-(n+1)} & \text{for type $A_n^{(1)}$,} \\
		o(r)(q-q^{-1})q^{-2(n-1)} & \text{for type $D_n^{(1)}$.}
	\end{cases}
\end{equation}
\end{thm}
\pf 
By definition \eqref{eq:borel action}, $\V$ is a weight vector and $\pfe_i\V=0$ for all $i\in I$. Let $M$ be the $U_q(\mf{b})$-submodule of $U_q^-(\varpi_r)_a$ generated by $\V$. By Propositions \ref{prop:highest l-weight of type A} (see Example \ref{ex:l-weight in A1} for type $A_1^{(1)}$) and \ref{prop:highest l-weight of type D}, $L_{r,ac_r}^+$ is a maximal quotient of $M$.
Note that the character of $U_q^-(\varpi_r)_a$ with respect to the action of $k_i$ for $i\in I$ is
\begin{equation}\label{eq:ch of prefundamental}
{\rm ch}\, U_q^-(\varpi_r)_a=\frac{1}{\prod_{\beta\in \Delta^+(w_r)}(1-e^{-\beta})},
\end{equation}
where $e^{-\beta}$ denotes an element of a basis of the group algebra of $\mathring{Q}$. 
On the other hand, the character of $L_{r,ac_r}^+$  is also equal to \eqref{eq:ch of prefundamental} by \cite[Theorem 7.5 and Lemma 7.6]{HJ} (cf. \cite{MY}), which implies that
\begin{equation*}
L_{r,ac_r}^+\cong M\cong U_q^-(\varpi_r)_a.
\end{equation*}
This completes the proof.
\qed 

\begin{rem} \label{rem:actions and limit construction}\
{\em 
\begin{enumerate}	
	\item In Section \ref{subsec:loop highest weight}, we have seen that
\begin{equation*}
	\E_{k\delta-\alpha_i} \V = 
	\begin{cases}
		q^{M} a (f_r)  & \text{if $k = 1$ and $i = r$,} \\
		0 & \text{otherwise,}
	\end{cases}
\end{equation*}
where $M$ is given in Corollaries \ref{cor:action of real root vector} and \ref{cor:action of real root vector for type D}.
Since the Drinfeld generator $x_{i, k}^-$ is given by 
\begin{equation*}
	x_{i, k}^- = -o(i)^k C^{-k}  k_i \E_{k\delta-\alpha_i}
\end{equation*}
for $k > 0$ (e.g. see \cite[Lemma 1.5]{BCP}), we see that the action of $x_{i, k}^-$ on $\V$ is $0$ for $k \ge 2$. By Theorem \ref{thm:main-2}, the description of the action of $x_{i, k}^-$ on $\V$ may be viewed as a special case of \cite[Proposition 7.3]{HJ} (cf. \cite[Theorem 6.3]{FH}).
	\vskip 2mm
	
	\item For an affine type, the limit construction of $L_{r, 1}^+$ by the family of KR modules $W^{(r)}_{s, q_i^{2s-1}}$ for $s \ge 1$ in \cite{HJ} is valid when $n_r = 1$ due to the convergence of $0$-action (see \cite[Proposition 7.4]{HJ}), where $q_i = q^{(\alpha_i, \alpha_i)/2}$ and $n_r$ is the multiplicity of $\alpha_r$ in the maximal root of $\mathring{\mf g}$.
In particular, for types $ADE$, $n_r = 1$ if and only if $\varpi_r$ is minuscule.
\end{enumerate}
}
\end{rem}
\vskip 1mm

\subsection{Realization of $L^-_{r,a}$} \label{subsec:realization L-}
Let us consider the realization of $L^-_{r,a}$.
Let $a \in \bR^\times$.
For $0\le i\le n$, we define the ${\bf k}$-linear operators $\pfe_i$ and $\pfk_i$ on $U_q^-({\mf g})$ by
\begin{equation} \label{eq:borel action-2}
\mathsf{e}_i(u)= 
\begin{cases}
e_i' (u) \quad &\text{ if } i\in I\,,  \\
a \xz \cdot u  \quad &\text{ if } i=0\,,
\end{cases}
\qquad
\mathsf{k}_i(u)=  
\begin{cases}
q^{(\al_i,\, \beta)}u \quad &\text{ if } i\in I\,,  \\
q^{ - (\theta,\, \beta)}u  \quad &\text{ if } i=0\,,
\end{cases} 
\end{equation}
for $u\in U_q^-({\mf g})_\beta$ ($\beta\in -Q_+$). 
Here $\xz$ is given in \eqref{eq:xz}.
The main difference from \eqref{eq:borel action} is that the action of $\mathsf{e}_0$ is given by the left multiplication of $\xz$ without $q^{-(\theta, \beta)}$.

\begin{ex} \label{ex:A1(1)-2}
{\em 
Let us consider the ${\bf k}$-linear operators $\pfe_i$ and $\pfk_i$ \eqref{eq:borel action-2} in type $A_1^{(1)}$.
We use the notation $S_{i,j}$ in Example \ref{ex:A1(1)}.
\vskip 2mm
First, we show that the $\pfe_i$ and $\pfk_j$ $(i, j \in \{ 0, 1 \})$ satisfy the defining relations of $U_q({\mf b})$.
In this example, we only verify $S_{0,1}(u) = 0$ for $u \in U_q^-({\mf g})_{\beta}$. 
It is almost identical to show that $S_{0,1}(u) = 0$, and it is straightforward to verify the relations for $\pfe_i$ and $\pfk_j$. 
We leave the details to the reader. 
\vskip 1mm

By definition \eqref{eq:borel action-2}, we get
\begin{equation*}
\begin{split}
	e_0^3 e_1 (u) &= a^3 f_1^3 e_1'(u)\,, \\
	e_0^2 e_1 e_0 (u) &= a^3 f_1^2 u + a^3 q^{-2} f_1^3 e_1'(u)\,, \\
	e_0 e_1 e_0^2 (u) &= a^3 (q^{-2} + 1) f_1^2 u + a^3 q^{-4} f_1^3 e_1'(u)\,, \\
	e_1 e_0^3 (u) &= a^3 ( q^{-4} + q^{-2} + 1) f_1^2 u + a^3 q^{-6} f_1^3 e_1'(u)\,.
\end{split}
\end{equation*}
Thus, we have $S_{0, 1}(u) = 0$.
\vskip 2mm

Next, let us compute the highest $\ell$-weight of $\V$ (cf. Example \ref{ex:l-weight in A1}).
By Lemma \ref{lem:inductive formula}, one can check that 
\begin{equation*}
	\E_{k\delta-\alpha_1} \V = (-1)^{k-1} q^{-2(k-1)} (q-q^{-1})^{k-1} a^k \cdot f,
\end{equation*}
where $k > 0$. Then it follows from \eqref{eq:Drinfeld generators in terms of affine root vectors} that the eigenvalue of $\psi_{1, k}^+$ on $\V$ is given by
\begin{equation*}
	\psi_{1, k}^+ \V = (-1)^k o(1)^k (q-q^{-1})^k q^{-2k} a^k \V.
\end{equation*}
Hence, the $\ell$-weight $\Psi(z)$ of $\V$ is given by
$$\left(\, \frac{1}{1+o(1)(q-q^{-1})q^{-2}z} \,\right).$$
\vskip 2mm

We have seen that $U_q^-({\mf g})$ is a $U_q({\mf b})$-module under \eqref{eq:borel action-2} and $L_{1, a}^-$ is a subquotient of $U_q^-({\mf g})$ up to a shift of spectral parameter.
}
\end{ex}

In this subsection, we will show that for a minuscule $\varpi_r$,\, $U_q^-(\varpi_r)$ is a $U_q(\mf{b})$-module under the actions in \eqref{eq:borel action-2} and it is isomorphic to $L_{r,1}^-$ up to a shift of spectral parameter.

\begin{thm} \label{thm: affine b-2}
The operators $\mathsf{e}_i$ and $\mathsf{k}_i$ for $0\le i\le n$ satisfy the defining relations of $U_q(\mf{b})$. This gives a representation $\rho_{r,a}^- : U_q(\mf{b}) \longrightarrow {\rm End}_{\bf k}(U_q^-(\mf{g}))$ defined by $\rho_{r,a}^-(e_i)=\mathsf{e}_i$ and $\rho_{r,a}^-(k_i)=\mathsf{k}_i$ for $0\le i\le n$. 
\end{thm}
\pf 
The proof is almost identical with the one of Theorems \ref{thm:main-1} and \ref{thm: affine b} except for verifying quantum Serre relations of $\pfe_0$ and $\pfe_i$ for $i \in J_1$.
Note that the case of type $A_1^{(1)}$ is done by Example \ref{ex:A1(1)-2}, so we consider the remaining cases.
\vskip 2mm

For $i \in J_1$ (see \eqref{eq:xz} below),
we claim that
\begin{align} 
	\pfe_0 \pfe_i^2 - (q+q^{-1})\pfe_i \pfe_0 \pfe_i + \pfe_i^2 \pfe_0 = 0,\label{eq:Serre relation 1'}\\
	\pfe_i \pfe_0^2 - (q+q^{-1})\pfe_0 \pfe_i \pfe_0 + \pfe_0^2 \pfe_i = 0.\label{eq:Serre relation 2'}
\end{align}

Let $u\in U_q^-(\mf{g})_\beta$ be given. Put $s = (\alpha_i, \beta)$ and $t = (\theta,\beta)$.
By Lemma \ref{lem:dual PBW properties 2}, we have
\begin{equation*}
\begin{split}
	\pfe_0\pfe_i^2(u) &= a \xz e_i'^2(u), \\
	\pfe_i \pfe_0 \pfe_i(u) &= a e_i'(\xz) e_i'(u) + aq^{-1} \xz e_i'^2(u), \\
	\pfe_i^2 \pfe_0(u) &= a(q+q^{-1})e_i'(\xz) e_i'(u) + aq^{-2}\xz e_i'^2(u),
\end{split}
\end{equation*}
which implies \eqref{eq:Serre relation 1'}. Similarly, we have 
\begin{equation*}
\begin{split}
	\pfe_i \pfe_0^2(u) &= a^2 q^{-1} (q+q^{-1})e_i'(\xz)\xz u + a^2 q^{-2}\xz^2 e_i'(u),\\
	\pfe_0 \pfe_i \pfe_0 (u) &= a^2  q^{-1} e_i'(\xz)\xz u + a^2q^{-1}\xz^2e_i'(u),\\
	\pfe_0^2 \pfe_i (u) &= a^2  \xz^2e_i'(u),
\end{split}
\end{equation*}
which implies \eqref{eq:Serre relation 2'}.
\qed

\begin{rem}
{\em 
For type $A_n^{(1)}$ with $r = 1$ or $n$, one can check that for $u\in U_q^-({\mf g})_\beta$, 
\begin{equation*}
	\pfe_r \pfe_0 (u) = q^{-1} \pfe_0 \pfe_r (u),
\end{equation*}
which also implies the relations \eqref{eq:Serre relation 1'} and \eqref{eq:Serre relation 2'} (cf. Remark \ref{rem:boundary in An(1)}).
}
\end{rem}
\vskip 2mm

In Section \ref{subsec:b-actions}, we have seen that $U_q^-(\varpi_r)$ is invariant under $e_i'$ for $i \in I$. Note that $U_q^-(\varpi_r)$ is also invariant under $\pfe_0$ in \eqref{eq:borel action-2}.
By abuse of notation, we also denote by $U_q^-(\varpi_r)_a$ the $U_q(\mf{b})$-module $U_q^-(\varpi_r)$ defined by \eqref{eq:borel action-2}.
\vskip 1mm

The second main result of this paper is as follows.
\begin{thm}\label{thm:main-3}
For a minuscule $\varpi_r$ and $a\in {\bf k}^\times$, we have 
\begin{equation*}
U_q^-(\varpi_r)_a \cong L_{r,\,-ac_r}^-,
\end{equation*} 
as a $U_q(\mf{b})$-module, 
where $c_{r}$ is given as in \eqref{eq:c_r}.
\end{thm}
\vskip 1mm

\begin{rem}
{\em 
It is shown in \cite{HJ} that $L_{r, a}^+$ is a $\sigma$-twisted dual of $L_{r, a}^-$ as a representation of the asymptotic algebra. 
Furthermore, the $q$-characters of $L_{r, a}^-$ and $L_{r, a}^+$ are far from being the same (e.g. see \cite[Remark 4.4]{FH}) while the characters of $L_{r, a}^-$ and $L_{r, a}^+$ with respect to $k_i\,(i \in I)$ are equal and independent of the choice of the parameter $a \in \bR^{\times}$ \cite[Theorem 6.4]{HJ}.
In this sense, Theorem \ref{thm:main-3} looks quite interesting. However,
it is not clear to us yet how the $U_q(\mf{b})$-actions \eqref{eq:borel action} and \eqref{eq:borel action-2} on $U_q^-(\varpi_r)$ result in these two different $U_q({\mf b})$-module structures.
}
\end{rem}

In the remainder of this section, we prove Theorem \ref{thm:main-3}. 
First, we compute the action of $\E_{k\delta-\alpha_i}$ $(i \in I,\, k > 0)$ on $\V$ as in Section \ref{subsec:loop highest weight}, where $\V$ is the weight vector in $U_q^-(\varpi_r)_a$ of weight $0$. The computation gives the $\ell$-weight of $\V$ by \eqref{eq:Drinfeld generators in terms of affine root vectors}, which enable us to apply the argument in the proof of Theorem \ref{thm:main-2}.

We remark that the computation for the $\ell$-weight of $\V$ in this case is more involved than in Section \ref{subsec:loop highest weight} since there is no cancellation as in the previous case when computing the action of $\E_{k\delta-\alpha_i}$ (see Example \ref{ex:A1(1)-2} for type $A_1^{(1)}$).
Let us illustrate this in the case of type $D_4^{(1)}$ with $r =4$ before considering a general case.

\begin{ex} \label{ex:l-weight in D4}
{\em 
Let us recall Example \ref{ex:root vector in D4}. We claim that 
\begin{equation*}
	\E_{k\delta-\alpha_i} \V = 
	\begin{cases}
		(-1)^{k-1} q^{-6k+2} (q-q^{-1})^{k-1} a^k (f_4) & \text{if $i = 4$,} \\
		0 & \text{otherwise.}
	\end{cases}
\end{equation*}
We proceed by induction on $k$.
Thanks to Lemma \ref{lem:dual PBW properties 1} and Lemma \ref{lem:inductive formula}, it suffice to verify the case of $i=4$ (cf. Example \ref{ex:l-weight in A3}).
\vskip 2mm

{\it Step 1.} 
By Lemma \ref{lem:dual PBW properties 1} and Example \ref{ex:root vector in D4}, we have
\begin{equation*}
\begin{split}
	\E_{\delta-\alpha_4} \V &= q^{-4} \pfe_2 \pfe_1 \pfe_3 \pfe_2 \pfe_0 \V = q^{-4}a (f_4), \\
	\E_{\delta-\alpha_4} (f_4) &= q^{-4} \pfe_2 \pfe_1 \pfe_3 \pfe_2 \pfe_0 (f_4) = q^{-4}a (f_4^2).
\end{split}
\end{equation*}
\vskip 1mm

{\it Step 2.}
Suppose that the claim holds for $k > 1$.
By Lemma \ref{lem:inductive formula} and the induction hypothesis together with {\it Step 1}, we have
\begin{equation*}
\begin{split}
	(q+q^{-1}) \E_{(k+1)\delta-\alpha_4} \V
	&=
	-C_k ( q^2 - q^{-2}  ) (f_4) = -C_k (q-q^{-1}) (q+q^{-1}) (f_4),
\end{split}
\end{equation*}
where $C_k = (-1)^{k-1} q^{-6(k+1)+2}(q-q^{-1})^{k-1}a^{k+1}.$
Hence, the claim is proved.
\vskip 1mm

By \eqref{eq:Drinfeld generators in terms of affine root vectors}, we have 
\begin{equation*}
	\psi_{i, k}^+ \V =
	\begin{cases}
		(-1)^k o(4)^k q^{-6k}(q-q^{-1})^k a^k \V  & \text{if $i = 4$,} \\
		0 & \text{otherwise,}
	\end{cases}
\end{equation*}
for $k > 0$.
Hence we obtain the $\ell$-weight $\Psi(z)$ of $\V$ as follows:
\begin{equation*}
	\Psi(z) = \left( 1,\, 1,\, 1,\, \frac{1}{1+ac_4 z} \right),
\end{equation*}
where $c_4 = o(4)q^{-6}(q-q^{-1})$.
}
\end{ex}

\subsubsection{\bf Type $A_n^{(1)}$}
Let us consider the case of type $A_n^{(1)}$ $(n \ge 2)$ with $r \in I$ (see Example \ref{ex:A1(1)-2} for type $A_1^{(1)}$).

\begin{lem} \label{lem:action of E on alpha for negative}
The action of $\E_{\delta-\alpha_i}$ on $f_r \in U_q^-(\varpi_r)_a$ is given by
\begin{equation*}
	\E_{\delta-\alpha_i} \left( f_r \right) 
	=
	\begin{cases}
		(-1)^{n-1} q^{-n+1} a  (f_r^2) & \text{if $i = r$,} \\
		0 & \text{if $i \neq r$.}
	\end{cases}
\end{equation*}
\end{lem}
\pf
By Lemma \ref{lem:formula of E} and \eqref{eq:borel action-2}, we have
\begin{equation*}
\begin{split}
	\E_{\delta-\alpha_r}\left( f_r \right)
	&=
	(-q^{-1})^{n-1} (\pfe_{r+1}\dots \pfe_{n-1}\pfe_n)(\pfe_{r-1}\dots \pfe_2\pfe_1)\pfe_0 \left( f_r \right) \\
	&=
	(-q^{-1})^{n-1}a (\pfe_{r+1}\dots \pfe_{n-1}\pfe_n)(\pfe_{r-1}\dots \pfe_2\pfe_1) \left( \xz f_r \right) \\
	&= 
	(-q^{-1})^{n-1}a (\pfe_{r+1}\dots \pfe_{n-1}\pfe_n) \left( F^{\rm up} (\ep_{r} - \ep_{n+1} ) f_r \right) \\
	&=
	(-q^{-1})^{n-1}a  \left(  f_r^2 \right).
\end{split}
\end{equation*}
The case of $i \neq r$ follows from Lemma \ref{lem:dual PBW properties 1} (see also Remark \ref{rem:minuscule case}).
\qed

\begin{lem} \label{lem:action of E on 1 for negative}
The action of $\E_{k\delta-\alpha_i}$ on $\V \in U_q^-(\varpi_r)_a$ $(k \ge 1)$ is given as follows:
\begin{equation*}
	\E_{k\delta-\alpha_i} \V =
	\begin{cases}
		\,\, (-1)^{kn-1} q^{-k(n+1)+2} (q-q^{-1})^{k-1} a^k ( f_r ) & \text{if $i = r$,} \\
		\,\, 0 & \text{if $i \neq r$.}
	\end{cases}
\end{equation*}
\end{lem}
\pf
The proof is by induction on $k$.
When $i \neq r$, 
by the same argument in the proof of Corollary \ref{cor:action of real root vector}, we have
$
	\E_{\delta-\alpha_i} \V = 0,
$
and then it follows from Lemma \ref{lem:inductive formula} that $\E_{k\delta-\alpha_i} \V = 0$ for $k \ge 1$.

Let us consider the case of $i = r$.
For $k = 1$, by Lemma \ref{lem:formula of E},
\begin{equation*}
\begin{split}
	\E_{\delta-\alpha_r} \V &= (-q^{-1})^{n-1} (\pfe_{r+1} \dots \pfe_{n-1} \pfe_n) (\pfe_{r-1} \dots \pfe_2 \pfe_1) \pfe_0 \V \\
	&= (-q^{-1})^{n-1} a (\pfe_{r+1} \dots \pfe_{n-1} \pfe_n) (\pfe_{r-1} \dots \pfe_2 \pfe_1) (\xz) \\
	&= (-q^{-1})^{n-1} a (\pfe_{r+1} \dots \pfe_{n-1} \pfe_n) ( F^{\rm up} (\ep_r - \ep_{n+1}) ) \\
	&= (-q^{-1})^{n-1} a (f_r).
\end{split}
\end{equation*}
Hence we have
\begin{equation} \label{eq:initial step for E}
\E_{\delta-\alpha_r} \V = (-1)^{n-1} q^{-n+1}a (f_r).
\end{equation}

Assume that the formula of $\E_{  k \delta-\alpha_r}\V$ holds for $k \ge 1$.
Let us compute the followings to proceed by induction on $k$.
By induction and \eqref{eq:initial step for E}, one can check that 
		\begin{equation*}
		\begin{split}
			\E_{\delta-\alpha_r} \E_{\alpha_r} \E_{k\delta-\alpha_r} \V 
			&= 
			C_k \cdot (f_r)
			= \E_{k\delta-\alpha_r} \E_{\alpha_r} \E_{\delta-\alpha_r} \V,
		\end{split}
		\end{equation*}
		where $C_k=(-1)^{(k+1)n-2} q^{-(k+1)(n+1)+4} a^{k+1} (q-q^{-1})^{k-1}$.
Also by induction and Lemma \ref{lem:action of E on alpha for negative}, one can check that
		\begin{equation*}
		\begin{split}
			\E_{\alpha_r} \E_{\delta-\alpha_r} \E_{k\delta-\alpha_r} \V
			&= C_k \cdot (1+q^{-2}) (f_r).
		\end{split}
		\end{equation*}
By Lemma \ref{lem:inductive formula} and the above computation, we get
\begin{equation*}
\begin{split}
	\E_{(k+1)\delta-\alpha_r}\V
	&= -C_k q^{-2}(q-q^{-1}) (f_r).
\end{split}
\end{equation*}
Hence, we have the desired formula of $\E_{(k+1)\delta-\alpha_r}\V$. This completes the induction.
\qed

\subsubsection{\bf Type $D_n^{(1)}$}
Let us assume that $n \ge 4$ and $r \in \{ 1, n-1, n\}$.

\begin{lem} \label{lem:action of E on alpha for negative in type D}
The action of $\E_{\delta-\alpha_i}$ on $f_r \in U_q^-(\varpi_r)_a$ is given by
\begin{equation*}
	\E_{\delta-\alpha_i} \left( f_r \right) 
	=
	\begin{cases}
		q^{-2n+4} a  ( f_r^2 ) & \text{if $i = r$,} \\
		0 & \text{if $i \neq r$.}
	\end{cases}
\end{equation*}
\end{lem}
\pf
Using Lemma \ref{lem:formula of E for type D}, the formula follows from the same argument as in the proof of Lemma \ref{lem:action of E on alpha for negative}.
\qed

\begin{lem} \label{lem:action of E on 1 for negative in type D}
The action of $\E_{k\delta-\alpha_i}$ on $\V \in U_q^-(\varpi_r)_a$ $(k \ge 1)$ is given as follows:
\begin{equation*}
	\E_{k\delta-\alpha_i} \V =
	\begin{cases}
		\,\, (-1)^{k-1} q^{-2k(n-1)+2} (q-q^{-1})^{k-1} a^k (f_r) & \text{if $i = r$,} \\
		\,\, 0 & \text{if $i \neq r$.}
	\end{cases}
\end{equation*}
\end{lem}
\pf
The proof is by induction on $k$ as in Lemma \ref{lem:action of E on 1 for negative}.
Then the necessary computation for induction is almost identical to the one in the proof of Lemma \ref{lem:action of E on 1 for negative} by using Lemma \ref{lem:formula of E for type D} and Lemma \ref{lem:action of E on alpha for negative in type D}. 
The details are left to the reader. 
\qed

\subsubsection{Proof of Theorem \ref{thm:main-3}} %\label{sec:pfof 2nd main theorem}
The character of $U_q^-(\varpi_r)_a$ is equal to \eqref{eq:ch of prefundamental} and it is also the one of $L_{r, a}^-$ for any $a \in {\bf k}^\times$ (cf. \cite[Theorem 6.4]{HJ}). By \eqref{eq:Drinfeld generators in terms of affine root vectors}, Lemma \ref{lem:action of E on 1 for negative} and Lemma \ref{lem:action of E on 1 for negative in type D},
we obtain the $\ell$-weight $\Psi=(\Psi_i(z))_{i\in I}$ of $\V$ as follows:
\begin{equation*}
	\Psi_i(z) = 
	\begin{cases}
		\,\, \displaystyle \frac{1}{1+ac_rz} & \text{if $i = r$,} \\
 		\,\, 1 & \text{if $i \neq r$.}
	\end{cases}
\end{equation*}
Then we may apply the argument in the proof of Theorem \ref{thm:main-2} to conclude that $U_q^-(\varpi_r)_a \cong L_{r,\,-ac_r}^-$ as a $U_q(\mf{b})$-module.
\qed

\begin{rem} \label{rem:characters}
{\em 
For types $A_n^{(1)}$ $(n \ge 1)$ and $D_n^{(1)}$ $(n \ge 4)$ with a minuscule $\varpi_r$, we have proved that $U_q^-(\varpi_r)$ has a $U_q({\mf b})$-module structure by which it is isomorphic to $L_{r, a}^{\pm}$ up to a spectral parameter shift.
It is natural to ask whether our approach could be extended to the non-symmetric types or a non-minuscule $\varpi_r$.
For types $B_n^{(1)}$, $C_n^{(1)}$, $D_n^{(1)}$ with a low rank $n$ and exceptional types $E_6^{(1)}$, $E_7^{(1)}$, $E_8^{(1)}$, $F_4^{(1)}$ and $G_2^{(1)}$ based on computer experiments, we observe that for any $r \in I$, the ordinary character of $U_q^-(\varpi_r)$ is given as follows:
\begin{equation} \label{eq:product formula}
	 \frac{1}{\prod_{\beta \in \Delta^+(w_r)} (1-e^{-\beta})^{[\beta]_r}},
\end{equation}
where $[\beta]_r \in \mathbb{Z}_+$ is the multiplicity of $\alpha_r$ in the sum $\beta = \sum_{s \in I} [\beta]_s \alpha_s$.
Indeed, the formula \eqref{eq:product formula} is equal to the one of the limit of the normalized characters of KR modules $W^{(r)}_{s,\, q_i^{-2s+1}}$ conjectured by Mukhin-Young \cite{MY} (cf. \cite{HJ}), which is proved recently in \cite{L19} for all untwisted types except for the case of type $E_8^{(1)}$ with $r \notin \{ 4, 8 \}$. Here $q_i = q^{(\alpha_i, \alpha_i)/2}$.
We expect that there exists a $U_q({\mf b})$-module structure on $U_q^-(\varpi_r)$ similar to \eqref{eq:borel action} (resp. \eqref{eq:borel action-2}) for all untwisted types and $r \in I$, which is isomorphic to $L_{r, a}^+$ (resp. $L_{r, a}^-$) up to a spectral parameter shift.
}
\end{rem}

\begin{rem}
{\em
Recall that in \cite{HJ}, Hernandez and Jimbo construct a module over an asymptotic algebra of $U_q({\mf g})$ by taking limits on the actions of Drinfeld generators of KR modules, which gives the prefundamental representation $L_{r, 1}^-$ for any $r \in I$. It is not yet clear to us how this limit construction is related to the realization of $L_{r, a}^\pm$ in this paper, where the actions of $e_i$ and $k_i$ $(i \in I)$ are defined more directly on the space $U_q^-(\varpi_r)$, and it would be interesting to clarify the connection between these two realizations.
}
\end{rem}

\section{Combinatorial realization of $L_{r, a}^+$}\label{sec:combinatorial realization}
In this section, as a byproduct of Theorem \ref{thm:main-2}, we give a combinatorial realization of $L_{r, a}^+$ in terms of the Lusztig data of the dual PBW vectors.

\begin{rem} \label{rem:combinatorial description of L-}
{\em 
We may try to give a combinatorial realization of $L_{r, a}^-$ from Theorem \ref{thm:main-3}. 
However, it doesn't seem to easy to describe the operator $\pfe_0$ in \eqref{eq:borel action-2} explicitly since the commutation relations associated to $\xz$ \eqref{eq:xz} may be complicated in general.
For this reason, we consider the combinatorial realization of $L_{r, a}^+$ only.
}
\end{rem}

\subsection{Type $A_n^{(1)}$}\label{subsec:type A} \
Let us identify the weight lattice $\mathring{P}$ for $\mathring{\mf{g}}$ with the abelian group generated by $\ep_i$ for $1\le i\le n+1$ subject to the relation $\ep_1+\cdots+\ep_{n+1}=0$ so that $\alpha_i=\ep_i-\ep_{i+1}$ for $1
\le i\le n$. Let $r \in I$ be given.
In this case, we have
$\Delta^+(w_r)=\{\,\ep_i-\ep_j\,|\,1\le i\le r <j \le n+1\,\}$.

Let
$$
\mc{C}(n, r)=\{\,{\bf c}=(c_{i,j})_{1\le i\le r <j \le n+1}\,|\, c_{i,j}\in\Z_+ \,\},
$$
and let 
$$
\mc{U}(n, r)=\bigoplus_{{\bf c} \, \in \, \mc{C}(n, r)}{\bf k} [{\bf c}]
$$
be the $\bR$-vector space with basis $\{\,[{\bf c}]\,|\,{\bf c}\in \mc{C}(n, r)\,\}$ parametrized by $\mc{C}(n, r)$.
We may identify ${\bf c}$ with an $r\times (n-r+1)$ matrix with non-negative integral entries, where the row index $i$ are decreasing from top to bottom. For example, when $n=6$ and $r=3$,
\begin{equation*}
{\bf c}=
\begin{pmatrix}
 c_{3,4}& c_{3,5} & c_{3,6} & c_{3,7} \\
 c_{2,4}& c_{2,5} & c_{2,6} & c_{2,7} \\
 c_{1,4}& c_{1,5} & c_{1,6} & c_{1,7} \\
\end{pmatrix}
\in 
\mc{C}(6, 3)\,.
\end{equation*}
By convention, if ${\bf c}$ is an $r\times(n-r+1)$ matrix with integral entries but ${\bf c} \notin \mc{C}(n, r)$, then we set $[{\bf c}]=0$ in $\mc{U}(n, r)$.
For $1\le k\le r <l \le n+1$, let ${\bf 1}_{k,l}=(c_{i,j}^{k,l})$ such that $c_{i,j}^{k,l}=\delta_{ik}\delta_{jl}$.

\begin{rem}
{\em 
Throughout this section, we use the notation $[m]_q$ $(m \in \mathbb{Z}_+)$ to distinguish the notation $[{\bf c}]$ for ${\bf c} \in \mathcal{C}(n, r)$.
}
\end{rem}

Fix $a\in {\bf k}^\times$.
For ${\bf c} = (c_{i,j}) \in \mc{C}(n, r)$, we set 
$$
{\rm wt}({\bf c}) 
= - \sum_{1\le i\le r <j \le n+1} c_{i,j} (\al_i + \cdots + \al_{j-1}) 
= - \sum_{1\le i\le r <j \le n+1} c_{i,j} (\ep_i - \ep_{j})\,,
$$
and define
{\allowdisplaybreaks
\begin{gather} \label{eq:type A actions}
e_i [{\bf c}] = 
\begin{cases} 
\,\,\, a q^{ \sum_{1 < k\le r} c_{k,n+1} + \sum_{r < l \le n+1} c_{1,l}}
[{\bf c}+{\bf 1}_{1,n+1}] &\text{ if } i=0,\vspace{0.2cm}\\
\,\, \displaystyle  \sum_{r < l \le n+1}  q^{ a_{i,l-1}({\bf c}) }  [c_{i,l}]_q [{\bf c}+{\bf 1}_{i+1,l}-{\bf 1}_{i,l}] &\text{ if } 1 \le i< r,\vspace{0.2cm}  \\
\,\,\, [c_{r,r+1}]_q [{\bf c}-{\bf 1}_{r,r+1}] &\text{ if } i=r,\vspace{0.2cm} \\
\,\, \displaystyle \sum_{ 1 \le k \le r } q^{ b_{k,i}({\bf c}) }  [c_{k,i+1}]_q
[{\bf c} + {\bf 1}_{k,i} - {\bf 1}_{k,i+1}]  &\text{ if } r < i \le n,  
\end{cases}
\end{gather}
}
{\allowdisplaybreaks
\begin{gather} \label{eq:type A actions k_i}
k_i [{\bf c}] = 
\begin{cases}
\,\, q^{(\al_i,\, {\rm wt}({\bf c}))} [{\bf c}] \quad &\text{ if } i\in I\,,  \\
\,\, q^{ - (\theta,\, {\rm wt}({\bf c}))}[{\bf c}] \quad &\text{ if } i=0\,,
\end{cases} 
\end{gather}
}
where we understand $c_{k,l}=0$ if $k >r $ or $l \le r$ and
\begin{align*}
a_{i,l}({\bf c}) &= \sum_{ r< t \le l } (c_{i+1,t}-c_{i,t}), \quad
b_{k,i}({\bf c}) = \sum_{ k< t \le r } (c_{t,i}-c_{t,i+1}).
\end{align*}
Here we assume that both $a_{i, r}({\bf c})$ and $b_{r, i}({\bf c})$ are $0$.
We denote by $\mc{U}(n, r)_a$ the space $\mc{U}(n, r)$ with the actions \eqref{eq:type A actions} and \eqref{eq:type A actions k_i} associated with $a \in {\bf k}^\times$.

\begin{lem} \label{lem:r th derivation and root vectors}
Let $\beta=\ep_k-\ep_l$ given with $1\le k\le r<l\le n+1$. 
For $i\in I$ and $c\ge 1$, we have
\begin{equation*}
\begin{split}
	e_i'\, F^{\rm up}(c\beta)
	&=
	\begin{cases}
			\,\, [c]_q\, F^{\rm up}((c-1)\beta) & \text{ if $i=r$ and $(k,l)=(r,r+1)$}, \\
			\,\, [c]_q\, F^{\rm up} \left(\beta-\alpha_i\right)F^{\rm up} ((c-1)\beta) & \text{ if $i\neq r$ and {\rm ($k=i$ or $l=i+1$)}},\\
			\,\, 0 & \text{otherwise},
	\end{cases}
\end{split}
\end{equation*}
where we understand $F^{\rm up}(0)= 1$.
\end{lem}
\pf 
We prove only the case when $i=r$ since the other cases are proved by similar computation as in \eqref{eq:computation e_r on F} below and Lemma \ref{lem:dual PBW properties 1}(2).

First, we consider the case of $\beta=\alpha_r$. We use induction on $c$.
The initial step follows from Lemma \ref{lem:dual PBW properties 1}(2).
Suppose $c > 1$. 
By \cite[Proposition 4.26]{Ki12}, one can verify that
\begin{equation} \label{eq:computation e_r on F}
\begin{split}
	e_r'\,F^{\rm up}(c\alpha_r) &= q^{c-1} e_r'\left( F^{\rm up} (\alpha_r) F^{\rm up} ((c-1)\alpha_r)\right) \\
						 &= q^{c-1} \big( F^{\rm up} ((c-1)\alpha_r) +  q^{-2} F^{\rm up} (\alpha_r) \left( [c-1]_q F^{\rm up} ( (c-2)\alpha_r ) \right) \big) \\
						 &= q^{c-1} \big( F^{\rm up} ((c-1)\alpha_r) +  q^{-1} [c-1]_q F^{\rm up} ( (c-1)\alpha_r ) \big) \\
						 &= [c]_q\, F^{\rm up} ((c-1)\alpha_r).
\end{split}
\end{equation}
Next, suppose that $\beta\neq \alpha_r$.
If $c = 1$, then we have $e_r' F^{\rm up}(\beta) = 0$ by Lemma \ref{lem:dual PBW properties 1}(2).
Since $e_r'$ is a derivation on $U_q^-(\mf{g})$, we have $e_r'\, F^{\rm up}(c\beta)=0$
by induction on $c$.
\qed\vskip 2mm

The following is a well-known identity which can be checked directly by using \eqref{eq:LS formula}.
\begin{lem} \label{lem:commuting root vectors}
For $t<s\le r< i\le n+1$ and $c\ge 1$, we have
\begin{equation*}
		\begin{split}
			\hspace{2cm} 
			F^{\rm up} (\ep_s - \ep_{i+1}) F^{\rm up} (\ep_t - \ep_i) &= F^{\rm up} (\ep_t - \ep_i) F^{\rm up} (\ep_s- \ep_{i+1}), \\
			F^{\rm up} (c(\ep_t - \ep_i)) F^{\rm up} (\ep_s - \ep_i) &= q^{-c} F^{\rm up} (\ep_s - \ep_i)  F^{\rm up} (c(\ep_t - \ep_i)). 
			\hspace{2cm} \qed
		\end{split}
		\end{equation*}	
\end{lem}
\vskip 1mm

\begin{thm}\label{thm:combinatorial realization for type A}
For $a\in {\bf k}^\times$ and $r \in I$, $\mc{U}(n, r)_a$ becomes a $U_q(\mf{b})$-module with respect to \eqref{eq:type A actions} and \eqref{eq:type A actions k_i}. 
Furthermore, we have
\begin{equation*}
	\mc{U}(n, r)_a \cong U_q^-(\varpi_r)_a \cong L_{r,ac_r}^+,
\end{equation*}
where $c_r$ is given as in \eqref{eq:c_r}.
\end{thm}
\pf
Let us take the reduced expression of $w_r$ in Remark \ref{rem:formula for w_i} whose convex order $\prec$ on $\Delta^+(w_r)$ is given by
\begin{equation} \label{eq:prec1}
\text{$\ep_i-\ep_j\prec \ep_k-\ep_l$ if and only if ($i>k$) or ($i=k$ and $j<l$)},
\end{equation}
for $\ep_i-\ep_j,\, \ep_k-\ep_l \in \Delta^+(w_r)$.
Let $\phi$ be the ${\bf k}$-linear isomorphism given by
\begin{equation}\label{eq:phi}
\xymatrixcolsep{3pc}\xymatrixrowsep{0pc}
\phi : \xymatrix{
   \mc{U}(n, r)_a \ \ar@{->}[r] & \quad \quad U_q^-(\varpi_r)_a \\
  [{\bf c}] \ar@{|->}[r] \quad & \,\,\, \displaystyle \prod_{\prec} F^{\rm up}\left(\,c_{i,j}(\ep_i - \ep_{j})\,\right)
},
\end{equation}
where $\prod_{\prec}$ is the ordered product with respect to $\prec$. 

We may take a reduced expression of $w_r$ such that the resulting convex order $\prec'$ on $\Delta^+(w_r)$ is given by
\begin{equation}\label{eq:prec2}
\text{$\ep_i-\ep_j\prec' \ep_k-\ep_l$ if and only if ($j<l$) or ($i>k$ and $j=l$)},
\end{equation}
for $\ep_i-\ep_j,\, \ep_k-\ep_l\in \Delta^+(w_r)$, and define a ${\bf k}$-linear isomorphism $\phi'$ as in \eqref{eq:phi} with respect to $\prec'$.
\vskip 2mm

We claim that $\phi=\phi'$.
Let us consider the following $r \times (n-r+1)$ matrix.
\begin{equation*}
\scalebox{0.9}{$\left( \begin{array}{cccc}
r & r+1 & \dots & n \\
r-1 & r & \dots & n-1 \\
\vdots & \vdots & \ddots & \vdots \\
1 & 2 & \dots & r
\end{array} \right)$}.
\end{equation*}
If we read the entries of the above matrix row by row (resp.~column by column) from top to bottom and then from left to right in each row (resp.~from left to right and then from top to bottom in each column), then the resulting sequence is associated to the reduced expression of $w_r$ whose convex order is $\prec$ \eqref{eq:prec1} (resp.~$\prec'$ \eqref{eq:prec2}).

One can check that these reduced expressions are equal up to $2$-braid moves (Recall that a $2$-braid move means $ij = ji$ for $i, j \in I$ such that $|i-j| > 1$). Hence, we have by Lemma \ref{lem:commuting root vectors},
\begin{equation*}
	\prod_{\prec} F^{\rm up}\left(\,c_{i,j}(\ep_i - \ep_{j})\,\right) = \prod_{\prec'} F^{\rm up}\left(\,c_{i,j}(\ep_i - \ep_{j})\,\right),
\end{equation*}
where $\ep_i-\ep_j \in \Delta^+(w_r)$. Hence, the claim is proved.
\vskip 2mm

Now, we show that $e_i$ and $k_i$ on $\mc{U}(n, r)_a$ coincide with $\mathsf{e}_i$ and $\mathsf{k}_i$ on $U_q^-(\varpi_r)$ under $\phi$.
It is clear that $\phi(k_i[{\bf c}]) = \pfk_i\,\phi([{\bf c}])$ for $0 \le i \le n$.
We also have $\phi(e_0[{\bf c}]) = \pfe_0\phi([{\bf c}])$ by definition and \cite[Proposition 4.26(2)]{Ki12}, and $\phi(e_r[{\bf c}]) = \pfe_r\phi([{\bf c}])$ by Lemma \ref{lem:r th derivation and root vectors}.
So it remains to check that $\phi(e_i[{\bf c}]) = \pfe_i\phi([{\bf c}])$ for the case when $i \neq 0, r$.
Let $[{\bf c}]\in \mc{U}(n, r)$ be given with ${\bf c}=(c_{i,j})\in \mc{C}(n, r)$.
\vskip 2mm

{\it Case 1.} Suppose that $1 \le i < r$. In this case, we use the ${\bf k}$-linear isomorphism $\phi'$ associated to $\prec'$ \eqref{eq:prec2}.
Let $y_0=\phi'([{\bf c}])$ and write
\begin{equation*}
\begin{split}
	y_0 = x_1 \cdot F^{\rm up} (c_{i+1, r+1}(\ep_{i+1} - \ep_{r+1})) F^{\rm up} (c_{i,r+1}(\ep_i - \ep_{r+1})) \cdot y_1,
\end{split}
\end{equation*}
where $x_1$ is the product of $F^{\rm up}(c_{u,v}(\ep_u-\ep_v))$ such that $\ep_u-\ep_v\prec' \ep_{i+1} - \ep_{r+1}$, and $y_1$ is the product of the other root vectors. 
Inductively, we write $y_{k}$ as follows:
\begin{equation*}
\begin{split}
	y_k &= x_{k+1} \cdot F^{\rm up} (c_{i+1, r+k+1} (\ep_{i+1} - \ep_{r+k+1})) F^{\rm up} (c_{i, r+k+1}(\ep_i - \ep_{r+k+1})) \cdot y_{k+1},
	\end{split}
\end{equation*}
where $x_{k+1}$ is the product of $F^{\rm up}(c_{u,v}(\ep_u-\ep_v))$ such that $\ep_u-\ep_v\prec' \ep_{i+1} - \ep_{r+1}$, and $y_{k+1}$ is the product of the other root vectors in $y_k$.
By Lemma \ref{lem:r th derivation and root vectors} and \cite[Proposition 4.26]{Ki12},

{\allowdisplaybreaks
\begin{equation} \label{eq:initial step-1}
\begin{split}
	e_i' y_0
	&= q^{c_{i+1, r+1}} x_1 \cdot F^{\rm up} (c_{i+1, r+1}(\ep_{i+1} - \ep_{r+1})) e_i' \Big( \, F^{\rm up} (c_{i,r+1}(\ep_i - \ep_{r+1})) \cdot y_1 \, \Big) \\
	&= q^{c_{i+1, r+1}} x_1 \cdot F^{\rm up} (c_{i+1, r+1}(\ep_{i+1} - \ep_{r+1})) \\
	&\quad \scalebox{1.5}{$\cdot$}\, \Big( [c_{i,r+1}]_q F^{\rm up} (\ep_{i+1} - \ep_{r+1} \big) F^{\rm up} ((c_{i,r+1}-1)(\ep_i - \ep_{r+1})) \cdot y_1 \\
	& \qquad \qquad \qquad \qquad \qquad \qquad \qquad \quad\,\,\, + q^{-c_{i,r+1}} F^{\rm up} (c_{i,r+1}(\ep_i - \ep_{r+1})) \cdot e_i'(y_1) \Big) \\
	&=[c_{i,r+1}]_q\,x_1 \cdot F^{\rm up}((c_{i+1, r+1}+1) (\ep_{i+1}-\ep_{r+1})) F^{\rm up} ((c_{i,r+1}-1)(\ep_i - \ep_{r+1})) \cdot y_1 \\
	& \quad + q^{c_{i+1, r+1}-c_{i, r+1}} x_1 \cdot F^{\rm up} (c_{i+1, r+1}(\ep_{i+1} - \ep_{r+1})) F^{\rm up} (c_{i,r+1}(\ep_i - \ep_{r+1})) \cdot e_i'(y_1)
\end{split}
\end{equation}}
Similarly, we have 
\begin{equation} \label{eq:general step-1}
\begin{split}
	e_i' y_k &=[c_{i,r+k+1}]_q\,x_{k+1} \cdot F^{\rm up}((c_{i+1, r+k+1}+1) (\ep_{i+1}-\ep_{r+k+1})) \\
	& \qquad \qquad \qquad \quad\,\,\,\,  \scalebox{1.5}{$\cdot$}\, F^{\rm up} ((c_{i,r+k+1}-1)(\ep_i - \ep_{r+k+1})) \cdot y_{k+1} \\
	& + q^{c_{i+1, r+k+1}-c_{i, r+k+1}} x_{k+1} \cdot F^{\rm up} (c_{i+1, r+k+1}(\ep_{i+1} - \ep_{r+k+1})) \\
	& \qquad \qquad \qquad \qquad \quad \,\,\,\,\,\,\,\, \scalebox{1.5}{$\cdot$}\, F^{\rm up} (c_{i,r+k+1}(\ep_i - \ep_{r+k+1})) \cdot e_i'(y_{k+1})
\end{split}
\end{equation}
Combining \eqref{eq:initial step-1} and \eqref{eq:general step-1}, we conclude that $\phi'(e_i [{\bf c}]) = \pfe_i\phi'([{\bf c}])$ for $1 \le i < r$.
\vskip 3mm

{\it Case 2.} Suppose that $r < i \le n$.
We may apply the same arguments as in {\it Case 1} by using $\phi$ or the linear order $\prec$ \eqref{eq:prec1}. We leave the details to the reader.
\qed

\begin{ex} \label{ex:combinatorics on actions in A}
{\em 
Let
\begin{equation*}
{\bf c} = 
\begin{pmatrix}
 c_{3,4}& c_{3,5} & c_{3,6} & c_{3,7} \\
 c_{2,4}& c_{2,5} & c_{2,6} & c_{2,7} \\
 c_{1,4}& c_{1,5} & c_{1,6} & c_{1,7} \\
\end{pmatrix} \in \mc{C}(6, 3).
\end{equation*}
In \eqref{eq:type A actions}, each coefficient of $[{\bf c}']$ appearing in the expansion of $e_i[{\bf c}]$ is determined by 
$c_{k, l}$'s which are relevant to $\alpha_i$. 
For example, when $i = 0$, it is determined by $c_{k, l}$'s except for the gray ones:
\begin{equation*}
\begin{pmatrix}
 \color{lightgray}{c_{3,4}}& \color{lightgray}{c_{3,5}} & \color{lightgray}{c_{3,6}} & c_{3,7} \\
 \color{lightgray}{c_{2,4}}& \color{lightgray}{c_{2,5}} & \color{lightgray}{c_{2,6}} & c_{2,7} \\
 c_{1,4}& c_{1,5} & c_{1,6} & c_{1,7} \\
\end{pmatrix}.
\end{equation*}
Here $(\theta, \epsilon_i-\epsilon_j) \neq 0$ for $i = 1$ or $j = 7$. Then we have 
\begin{equation*}
	e_0 [{\bf c}] =
	aq^{\sum_{1 < i \le 3} c_{i, 7} + \sum_{3 < j \le 7} c_{1, j}}
	\left[ \begin{pmatrix}
 c_{3,4}& c_{3,5} & c_{3,6} & c_{3,7} \\
 c_{2,4}& c_{2,5} & c_{2,6} & c_{2,7} \\
 c_{1,4}& c_{1,5} & c_{1,6} & c_{1,7} + 1\\
\end{pmatrix} \right].
\end{equation*}

As another example, when $i = 2$, the coefficients appearing in $e_i[{\bf c}]$ are determined by $c_{k, l}$'s except for the gray ones:
\begin{equation*}
\begin{pmatrix}
 c_{3,4}& c_{3,5} & c_{3,6} & c_{3,7} \\
 c_{2,4}& c_{2,5} & c_{2,6} & c_{2,7} \\
 \color{lightgray}{c_{1,4}}& \color{lightgray}{c_{1,5}} & \color{lightgray}{c_{1,6}} & \color{lightgray}{c_{1,7}} \\
\end{pmatrix}.
\end{equation*}
Then we have
\begin{equation*}
\begin{split}
	e_2 \left[ {\bf c} \right] &= 
	\quad \quad \quad
	[c_{2, 4}]_q
	\scalebox{0.7}{$
	\left[ \begin{pmatrix}
 c_{3,4}+1& c_{3,5} & c_{3,6} & c_{3,7} \\
 c_{2,4}-1& c_{2,5} & c_{2,6} & c_{2,7} \\
 c_{1,4}& c_{1,5} & c_{1,6} & c_{1,7} \\
\end{pmatrix} \right]
$}
+ 
q^{a_{2,4}({\bf c})} [c_{2, 5}]_q
\scalebox{0.7}{$
\left[  \begin{pmatrix}
 c_{3,4}& c_{3,5}+1 & c_{3,6} & c_{3,7} \\
 c_{2,4}& c_{2,5}-1 & c_{2,6} & c_{2,7} \\
 c_{1,4}& c_{1,5} & c_{1,6} & c_{1,7} \\
\end{pmatrix} \right]
$} \\ 
& \,+
q^{a_{2,5}({\bf c})} [c_{2, 6}]_q
	\scalebox{0.7}{$
	\left[  \begin{pmatrix}
 c_{3,4}& c_{3,5} & c_{3,6}+1 & c_{3,7} \\
 c_{2,4}& c_{2,5} & c_{2,6}-1 & c_{2,7} \\
 c_{1,4}& c_{1,5} & c_{1,6} & c_{1,7} \\
\end{pmatrix} \right]
$}
+
q^{a_{2,6}({\bf c})} [c_{2, 7}]_q
	\scalebox{0.7}{$
	\left[  \begin{pmatrix}
 c_{3,4}& c_{3,5} & c_{3,6} & c_{3,7} +1 \\
 c_{2,4}& c_{2,5} & c_{2,6} & c_{2,7} -1 \\
 c_{1,4}& c_{1,5} & c_{1,6} & c_{1,7} \\
\end{pmatrix} \right]
$},
\end{split}
\end{equation*}
where $\displaystyle a_{2,l}({\bf c}) = \sum_{3 < t \le l} (c_{3, t}-c_{2, t})$ and $a_{2, 3}({\bf c}) = 0$ by definition.
}
\end{ex}

\subsection{Type $D_n^{(1)}$}\label{subsec:type D} 
Let us consider the case of type $D_n^{(1)}$. Here we assume that $r = 1$ or $n$. Note that the case of $r = n-1$ is almost identical to the case of $r=n$.
We assume that the weight lattice of $\mathring{\mf{g}}$ is $\mathring{P} = \bigoplus_{i=1}^n \mathbb{Z} \epsilon_i$ with a symmetric bilinear form $(\,,\,)$ such that $(\epsilon_i, \epsilon_j) = \delta_{i,j}$ for $1 \le i,\,j \le n$, so that $\alpha_i = \epsilon_i - \epsilon_{i+1}$ $(1 \le i \le n-1)$ and $\alpha_n = \epsilon_{n-1} + \epsilon_n$.
Then we have
\begin{equation*}
	\Delta^+(w_r) =
	\begin{cases}
		\left\{ \, 
			\epsilon_1 - \epsilon_i \, | \, 1 < i \le n  
		\, \right\} 
		\cup 
		\left\{ \, 
		\epsilon_1 + \epsilon_i \, | \, 1 < i \le n 
		\, \right\} & \text{if $r = 1$}, \\
		\left\{ \,
			\epsilon_i + \epsilon_j \, | \, 1 \le i < j \le n
		\, \right\} & \text{if $r = n$}.
	\end{cases}
\end{equation*}

For each $\beta \in \Delta^+(w_r)$, we will denote the multiplicity of $\beta$ in a PBW vector by 
\begin{equation*}
		\begin{cases}
			c_{i, j} & \text{if $\beta = \epsilon_i - \epsilon_j$ for $1 \le i < j \le n$,} \\
			c_{\ov{j}, \ov{i}} & \text{if $\beta = \epsilon_i + \epsilon_j$ for $1 \le i < j \le n$,}
		\end{cases}
\end{equation*}
where we assume that $\ov{j} < \ov{i}$ for $1 \le i < j \le n$.

Set
\begin{equation*}
	\mc{C}(n, r) =
	\begin{cases}
		\left\{ \, {\bf c} =  \left( (c_{1, i})_{1 < i \le n},\, (c_{\ov{i}, \ov{1}})_{1 < i \le n} \right) \, \Big| \, c_{1, i},\, c_{\ov{i}, \ov{1}} \in \mathbb{Z}_+   \, \right\} & \text{if $r = 1$}, \\
		\left\{ \, {\bf c} = \left( c_{\ov{j}, \ov{i}} \right)_{1 \le i < j \le n} \, \Big| \, c_{\ov{j}, \ov{i}} \in \mathbb{Z}_+ \, \right\} & \text{if $r = n$,}
	\end{cases}
\end{equation*}
and let
$$
\mc{U}(n, r)=\bigoplus_{{\bf c} \, \in \, \mc{C}(n, r)}{\bf k} [{\bf c}]
$$
be the $\bR$-space with basis $\{\,[{\bf c}]\,|\,{\bf c}\in \mc{C}(n, r)\,\}$ parametrized by $\mc{C}(n, r)$, where we assume that $[{\bf c}] = 0$ if ${\bf c} \notin \mc{C}(n, r)$. 
We may identify ${\bf c}$ with a strictly upper triangular $n \times n$-matrix (resp. a pair of $1 \times (n-1)$-matrices) if $r = n$ (resp. $r = 1$).
For example, when $n = 4$ and $r = 4$,
\begin{equation*}
	{\bf c} = 
	\left( 
		\begin{array}{cccc}
			 & c_{\ov{4}, \ov{3}} & c_{\ov{4}, \ov{2}} & c_{\ov{4}, \ov{1}} \\
			 &                    & c_{\ov{3}, \ov{2}} & c_{\ov{3}, \ov{1}} \\
			 &                    &                    & c_{\ov{2}, \ov{1}} \\
			 &                    &                    &           
		\end{array}
	\right).
\end{equation*}

Fix $a\in {\bf k}^\times$. Let us define the weight function and action of $e_i$ on $\mc{U}(r)$ as follows: 
\vskip 3mm

{\it Case 1}. $r = 1$.
Put ${\bf 1}_{k, l} = \left( \left( c_{i,j}^{k,l} \right),\, \left( c_{\ov{j},\ov{i}}^{k,l} \right)\right)$ with $c_{a,b}^{c,d} = \delta_{ac} \delta_{bd}$.
For ${\bf c} = \left( (c_{1, i}),\, (c_{\ov{i}, \ov{1}}) \right) \in \mc{C}(n, 1)$, we set 
\begin{equation*}
	{\rm wt}({\bf c}) =
	-\sum_{1 < i \le n} \left( c_{1, i} (\epsilon_1 - \epsilon_i) + c_{\ov{i}, \ov{1}} (\epsilon_1 + \epsilon_i) \right)\,,
\end{equation*}	
and define
\begin{equation} \label{eq:type D actions with r=1}
\begin{split}
	e_i[{\bf c}] &=
	\begin{cases}
		 aq^{\sum_{2 < k \le n}( c_{1, k} + c_{\ov{k}, \ov{1}} )+c_{\ov{2}, \ov{1}}} [ {\bf c} + {\bf 1}_{\ov{2}, \ov{1}} ] & \text{if $i = 0$,} \\
		 [c_{1, 2}]_q[{\bf c}-{\bf 1}_{1, 2}] & \text{if $i = 1$,} \\
		 [c_{1, i+1}]_q[{\bf c}+{\bf 1}_{1, i}-{\bf 1}_{1,i+1}] + q^{c_{1, i}-c_{1,i+1}}[c_{\ov{i}, \ov{1}}]_q [{\bf c}+{\bf 1}_{\ov{i+1}, \ov{1}}-{\bf 1}_{\ov{i},\ov{1}}] & \text{if $1 < i \le n-1$,} \\
		 %[c_{1, n}]_q [{\bf c}+{\bf 1}_{1, n-1}-{\bf 1}_{1, n}] + q^{c_{1, n-1} - c_{1, n}} [c_{\ov{n-1}, \ov{1}}]_q [{\bf c} + {\bf 1}_{\ov{n}, \ov{1}} - {\bf 1}_{\ov{n-1}, \ov{1}}] & \text{if $i = n-1$,} \\
		 [c_{\ov{n}, \ov{1}}]_q [ {\bf c} + {\bf 1}_{1, n-1} - {\bf 1}_{\ov{n}, \ov{1}} ] + q^{c_{1, n-1} - c_{\ov{n}, \ov{1}}}[c_{\ov{n-1}, \ov{1}}]_q [{\bf c} + {\bf 1}_{1, n} - {\bf 1}_{\ov{n-1}, \ov{1}}] & \text{if $i = n$.} \\
	\end{cases}
\end{split}
\end{equation}
\vskip 3mm

{\it Case 2}. $r = n$.
Put ${\bf 1}_{k, l} = \left( c_{\ov{j}, \ov{i}}^{k, l} \right)$ with $c_{a,b}^{c,d} = \delta_{ac} \delta_{bd}$.
For ${\bf c} = (c_{\ov{j}, \ov{i}}) \in \mc{C}(n, n)$, we set 
\begin{equation*}
	{\rm wt}({\bf c}) =
	-\sum_{1 \le i < j \le n} c_{\ov{j}, \ov{i}} (\epsilon_i + \epsilon_j)\,,
\end{equation*}	
and define
\begin{equation} \label{eq:type D actions with r=n}
\begin{split}
	e_i[{\bf c}] &=
	\begin{cases}
		 aq^{\sum_{1 < j \le n}  c_{\ov{j}, \ov{1}} + \sum_{2 < j \le n} c_{\ov{j}, \ov{2}}} [{\bf c} + {\bf 1}_{\ov{2}, \ov{1}}]& \text{if $i = 0$,} \\
		\displaystyle\sum_{1 \le k < l \le n} q^{a_{k,l}^i({\bf c})} [c_{\ov{l}, \ov{k}}]_q 
		 \left(
		 	\delta_{k\, i+1}[{\bf c} + {\bf 1}_{\ov{l}, \ov{k}} - {\bf 1}_{\ov{l}, \ov{k-1}}]  + \delta_{li}[{\bf c} + {\bf 1}_{\ov{l+1},  \ov{k}}  - {\bf 1}_{\ov{l}, \ov{k}}] 
		 \right) & \text{if $1 \le i < n$,} \\
		  [c_{\ov{n-1}, \ov{n}}]_q [{\bf c} - {\bf 1}_{\ov{n}, \ov{n-1}}] & \text{if $i = n$,}
	\end{cases}
\end{split}
\end{equation}
where $a_{k, l}^i({\bf c})$ is given by
\begin{equation*}
	a_{k, l}^i({\bf c}) = 
		\begin{cases}
			\,\,\displaystyle\sum_{l < p} (c_{\ov{p}, \ov{i+1}} - c_{\ov{p}, \ov{i}}) & \text{if $k=i$,} \\
			\,\,\displaystyle\sum_{p=i+2}^{n} (c_{\ov{p}, \ov{i+1}} - c_{\ov{p}, \ov{i}}) + \sum_{q < l} ( c_{\ov{i+1}, \ov{q}} - c_{\ov{i}, \ov{q}} ) & \text{if $l=i$.}
		\end{cases}
\end{equation*}

Finally, we define the action of $k_i$ on $\mc{U}(n, r)$ by
\begin{equation} \label{eq:action of k}
	k_i[{\bf c}] =
	\begin{cases}
		 q^{(\al_i,\, {\rm wt}({\bf c}))} [{\bf c}] \quad &\text{ if } i\in I,  \\
		 q^{ - (\theta,\, {\rm wt}({\bf c}))}[{\bf c}] \quad &\text{ if } i=0.
	\end{cases}
\end{equation}

\begin{ex} \label{ex:combinatorics on actions in D}
{\em 
Let us consider the case of type $D_4^{(1)}$ with $r =1, 4$ and $i = 2$.
\vskip 1mm

(1) Let $r = 1$. For ${\bf c} = \left( (c_{1,2}, c_{1, 3}, c_{1, 4}), (c_{\ov{4}, \ov{1}}, c_{\ov{3}, \ov{1}}, c_{\ov{2}, \ov{1}}) \right) \in \mc{U}(4, 1)$, 
the action of $e_2$ is given by
\begin{equation*}
\begin{split}
	e_2 [{\bf c}] 
	&=[c_{1, 3}]_q \left[ \left( (c_{1,2}+1, c_{1, 3}-1, c_{1, 4}), (c_{\ov{4}, \ov{1}}, c_{\ov{3}, \ov{1}}, c_{\ov{2}, \ov{1}}) \right) \right] \\
	&\quad +\, q^{c_{1, 2}-c_{1, 3}} [c_{\ov{2}, \ov{1}}]_q \left[ \left( (c_{1,2}, c_{1, 3}, c_{1, 4}), (c_{\ov{4}, \ov{1}}, c_{\ov{3}, \ov{1}}+1, c_{\ov{2}, \ov{1}}-1) \right) \right].
\end{split}
\end{equation*}
\vskip 2mm

(2) Let $r = 4$ and 
\begin{equation*}
	{\bf c} = 
	\left( 
		\begin{array}{cccc}
			 & c_{\ov{4}, \ov{3}} & c_{\ov{4}, \ov{2}} & c_{\ov{4}, \ov{1}} \\
			 &                    & c_{\ov{3}, \ov{2}} & c_{\ov{3}, \ov{1}} \\
			 &                    &                    & c_{\ov{2}, \ov{1}} \\
			 &                    &                    &           
		\end{array}
	\right) \in \mc{U}(4, 4).
\end{equation*}
In this case, the action of $e_2$ is determined by $c_{\ov{j}, \ov{i}}$'s except for the gray ones:
\begin{equation*}
	\left( 
		\begin{array}{cccc}
			 & c_{\ov{4}, \ov{3}} & c_{\ov{4}, \ov{2}} & \color{lightgray}{c_{\ov{4}, \ov{1}}} \\
			 &                    & \color{lightgray}{c_{\ov{3}, \ov{2}}} & c_{\ov{3}, \ov{1}} \\
			 &                    &                    & c_{\ov{2}, \ov{1}} \\
			 &                    &                    &           
		\end{array}
	\right).
\end{equation*}
We have 
\begin{equation*}
	a_{2, 4}^2 ({\bf c}) = 0, \qquad a_{1, 2}^2 ({\bf c}) = c_{\ov{4}, \ov{3}} - c_{\ov{4}, \ov{2}}.
\end{equation*}
Hence, the action $e_2[{\bf c}]$ is given by 
\begin{equation*}
	e_2 [{\bf c}] = [c_{\ov{4}, \ov{2}}]_q 
	\scalebox{0.8}{$
	\left[ \left( 
		\begin{array}{cccc}
			 & c_{\ov{4}, \ov{3}} +1 & c_{\ov{4}, \ov{2}} - 1 & c_{\ov{4}, \ov{1}} \\
			 &                    & c_{\ov{3}, \ov{2}} & c_{\ov{3}, \ov{1}} \\
			 &                    &                    & c_{\ov{2}, \ov{1}} \\
			 &                    &                    &           
		\end{array}
	\right) \right] $}
	+
	q^{c_{\ov{4}, \ov{3}} - c_{\ov{4}, \ov{2}}} [c_{\ov{2}, \ov{1}}]_q
	\scalebox{0.8}{$
	\left[ \left( 
		\begin{array}{cccc}
			 & c_{\ov{4}, \ov{3}} & c_{\ov{4}, \ov{2}} & c_{\ov{4}, \ov{1}} \\
			 &                    & c_{\ov{3}, \ov{2}} & c_{\ov{3}, \ov{1}}+1 \\
			 &                    &                    & c_{\ov{2}, \ov{1}}-1 \\
			 &                    &                    &           
		\end{array}
	\right) \right] $}
\end{equation*}
}
\end{ex}
\vskip 1mm

We denote by $\mc{U}(n, r)_a$ the space $\mc{U}(n, r)$ with the above actions of $e_i$ and $k_i$ associated with $a \in {\bf k}^\times$.

\begin{thm}\label{thm:combinatorial realization for type D}
For $a\in {\bf k}^\times$ and $r \in \{ \, 1, n \, \}$, $\mc{U}(n, r)_a$ becomes a $U_q(\mf{b})$-module with respect to \eqref{eq:type D actions with r=1}, \eqref{eq:type D actions with r=n}, and \eqref{eq:action of k}. Furthermore, we have
\begin{equation*}
	\mc{U}(n, r)_a \cong U_q^-(\varpi_r)_a \cong L_{r,ac_r}^+,
\end{equation*}
where $c_r$ is given as in \eqref{eq:c_r}. 
\end{thm}
\pf
The proof is almost identical to the one of Theorem \ref{thm:combinatorial realization for type A}.
Let us explain briefly the proof of Theorem \ref{thm:combinatorial realization for type D}.
First, we define a ${\bf k}$-linear map from $\mc{U}(n, r)_a$ to $U_q^-(\varpi_r)_a$ as in \eqref{eq:phi}.
\vskip 2mm

{\it Case 1}. $r = 1$. 
Let us take the reduced expression of $w_1$ in Remark \ref{rem:formula for w_i}. In this case, the convex order $\prec$ of $w_1$ is given by
\begin{equation*}
\begin{split}
 \ep_1  - \ep_i & \prec \ep_1 + \ep_j, \\
\ep_1 - \ep_j \prec \ep_1-& \ep_k \quad \Longleftrightarrow \quad 
\text{$j < k$,} \\
\ep_1 + \ep_j \prec \ep_1 \,+&\, \ep_k \quad \Longleftrightarrow \quad 
\text{$j > k$,}
\end{split}
\end{equation*}
for $i, j, k \in I$. Then we define the ${\bf k}$-linear map $\phi_1$ from $\mc{U}(n, 1)_a$ to $U_q^-(\varpi_1)_a$ by using the above convex order $\prec$ as in \eqref{eq:phi}.
\vskip 2mm

{\it Case 2}. $r = n$. 
Let us take the reduced expression of $w_n$ in Remark \ref{rem:formula for w_i}. 
Then the convex order $\prec$ of $w_n$ is given by 
\begin{equation} \label{eq:prec1_D}
\begin{split}
\ep_i + \ep_j \prec \ep_k \,+&\, \ep_l \quad \Longleftrightarrow \quad 
\text{$(j>l)$ or $(j=l$, $i>k)$}
\end{split}
\end{equation}
for $1\leq i<j\leq n$ and $1\leq k<l\leq n$. 
Also, we may take a reduced expression of $w_n$ so that the corresponding convex order denoted by $\prec'$ is given by
\begin{equation} \label{eq:prec2_D}
\begin{split}
\ep_i+\ep_j \prec' \ep_k \,+&\, \ep_l \quad \Longleftrightarrow \quad 
\text{$(i>k)$ or $(i=k$, $j>l)$}
\end{split}
\end{equation}
for $1\leq i<j\leq n$ and $1\leq k<l\leq n$.
We define the ${\bf k}$-linear maps $\phi_n$ and $\phi_n'$ from $\mc{U}(n, n)_a$ to $U_q^-(\varpi_n)_a$ by using the above convex orders \eqref{eq:prec1_D} and \eqref{eq:prec2_D} as in \eqref{eq:phi}, respectively.
\vskip 2mm

We show that $\phi_n = \phi_n'$ as in type $A_n^{(1)}$.
Let us consider a upper triangular $(n-1) \times (n-1)$ matrix ${\bf D}_n$ such that its diagonal is given by $(n,\, n-1,\, n,\, n-1,\, \cdots)$
and the remaining non-zero part of ${\bf D}_n$ is of the form
\begin{equation*}
\scalebox{0.9}{$\left( \begin{array}{cccccc}
\asterisk& n-2 & n-3 & n-4 & \dots & 1 \\
 & \asterisk  & n-2 & n-3 & \dots & 2 \\
 & & \ddots & \ddots &\ddots & \vdots \\
 & &   & \asterisk &n-2 & n-3 \\
 & &  &   & \asterisk & n-2 \\
 & & & & & \asterisk 
\end{array} \right)$},\,
\end{equation*}
where $\asterisk$'s denote the diagonal entries.
For example, for type $D_5^{(1)}$, the matrix ${\bf D}_5$ is given by
\begin{equation*}
\scalebox{0.9}{$\left(
	\begin{array}{cccc}
		5 & 3 & 2 & 1 \\
		  & 4 & 3 & 2 \\
		  & & 5 & 3 \\
		  & & &  4\\
	\end{array}
	\right)$}.
\end{equation*}

If we read the entries of ${\bf D}_n$ row by row (resp.~column by column) from top to bottom and then from left to right in each row (resp.~from left to right and then from top to bottom in each column), then the resulting sequence is associated to the reduced expression of $w_n$ whose convex order is $\prec$ \eqref{eq:prec1_D} (resp.~$\prec'$ \eqref{eq:prec2_D}).
It is straightforward to check that these reduced expressions are equal up to $2$-braid moves. Thus, by Lemma \ref{lem:commuting root vectors}, we prove the claim.
\vskip 2mm

Now, we apply the argument in the proof of Theorem \ref{thm:combinatorial realization for type A} to the ${\bf k}$-linear maps $\phi_r$ $(r = 1,\,n)$.
In particular, we should remark that in the case of $r = n$, the computation to show that $\phi_n(e_k[{\bf c}]) = \pfe_k\phi_n([{\bf c}])$ for $k \neq 0,\, n$ is more involved.
Let us explain this case in more details.
\vskip 1mm

Let $r = n$ and $N = n^2 - n$.
For ${\bf c} = (c_{\ov{j}, \ov{i}})_{1 \le i < j \le n} \in \mathbb{Z}_+^N$, put
\begin{equation*}
	F({\bf c}) = \prod_{\prec} F^{\rm up} (c_{\ov{j}, \ov{i}} (\ep_i + \ep_j)) \in U_q^-(\varpi_n). 
\end{equation*}
Suppose that $k \in I \,\setminus\, \{ n \}$.
Then $F({\bf c})$ is rewritten as follows:
\begin{enumerate}
	\item For $k+1 < j \leq n$, we define
		\begin{equation*}
			x_j = y_j \cdot \left( \prod_{\substack{\prec \\ \text{$i = k$ or $k+1$}}} F^{\rm up} (c_{\ov{j},\, \ov{i}} (\ep_i + \ep_j)) \right) \cdot z_j,
		\end{equation*}
		where $y_j$ and $z_j$ are given by 
		\begin{equation*}
			 y_j = \prod_{\substack{\prec \\ \text{$i < k+1$}}} F^{\rm up} (c_{\ov{j},\, \ov{i}} (\ep_i + \ep_j)), \quad\,\,
			 z_j = \prod_{\substack{\prec \\ \text{$i > k$}}} F^{\rm up} (c_{\ov{j},\, \ov{i}} (\ep_i + \ep_j)).
		\end{equation*}
		Here $y_j$ and $z_j$ are assumed to be $1$, if they are not defined.
		\vskip 1mm
		
	\item For $1 \le i \le k+1$, we define
		\begin{equation*}
			x_i' = y_i' \cdot \left( \prod_{\substack{\prec \\ \text{$j = k$ or $k+1$}}} F^{\rm up} (c_{\ov{j},\, \ov{i}} (\ep_i + \ep_j)) \right),
		\end{equation*}
		where $y_i'$ is given by
		\begin{equation*}
			 y_i' = \prod_{\substack{\prec \\ \text{$j \neq k,\,k+1$}}} F^{\rm up} (c_{\ov{j},\, \ov{i}} (\ep_i + \ep_j)).
		\end{equation*}
		Here if $y_i'$ and the term in $x_i'$ except for $y_i'$  are not defined, then we assume that they are equal to $1$.
\end{enumerate}

From (1) and (2), we rewrite $F({\bf c})$ by
\begin{equation*}
\begin{split}
	F({\bf c}) &= (x_n \cdot x_{n-1} \cdot\, \cdots\, \cdot x_{k+2}) \cdot (x_{k+1}' x_k' \cdot\, \cdots \, \cdot x_1') \\
	&= (x_n \cdot x_{n-1} \cdot\, \cdots\, \cdot x_{k+2}) \cdot (x_{k+1}'' x_k'' \cdot\, \cdots \, \cdot x_1''),
\end{split}
\end{equation*}
where $x_i''$ is obtained from $x_i'$ by rewriting with respect to the convex order $\prec'$ \eqref{eq:prec2_D} or using the map $\phi_n'$. Then we have
\begin{equation*}
\begin{split}
	\pfe_k F({\bf c}) &= e_k' ( F({\bf c}) ) \\
	&= e_i' ( x_n \cdot x_{n-1} \cdot\, \cdots\, \cdot x_{k+2} ) \cdot x_{k+1}'' x_k'' \cdot\, \cdots \, \cdot x_1'' \\
	& \qquad + q^{\sum_{p=k+2}^{n} (c_{\ov{p}, \ov{k+1}} - c_{\ov{p}, \ov{k}})} ( x_n \cdot x_{n-1} \cdot\, \cdots\, \cdot x_{k+2} ) \cdot e_i'( x_{k+1}'' x_k'' \cdot\, \cdots \, \cdot x_1'' ).
\end{split}
\end{equation*}

Finally, we compute the above terms containing $e_i'$ as in {\it Case 1} in the proof of Theorem \ref{thm:combinatorial realization for type A}, and one can check that $\phi_n(e_k[{\bf c}]) = \pfe_k\phi_n([{\bf c}])$ for $k \in I \, \setminus \, \{ n \}$.
We leave the details to the reader.
\qed

\end{document}